\documentclass[12pt, reqno]{amsart}
\usepackage{amsthm}
\usepackage{amssymb}
\usepackage{amsmath}

\def\1{\hbox{1\kern-.35em\hbox{1}}}

\def\bigl#1{\left#1\rb{1.2pt}{\mbox{${^{^{^{\,\!\!\!\!}}}}$}}\right.}
\def\bigr#1{\left. \rb{1.2pt}{\mbox{${^{^{^{\,\!\!\!\!}}}}$}}\right#1}
\newtheorem{theorem}{Theorem}[section]
\newtheorem*{theorem*}{Theorem}
\newtheorem{lemma}[theorem]{Lemma}
\newtheorem{proposition}[theorem]{Proposition}
\newtheorem*{proposition*}{Proposition}

\newtheorem{definition}[theorem]{Definition}
\newtheorem{notation}[theorem]{Definition}\def\Notation{Definition\ }
\newtheorem{coro}[theorem]{Corollary}
\newtheorem{remark}[theorem]{Remark}

\newtheorem{example}[theorem]{Example}

\newtheorem{desc}[theorem]{Description}
\newtheorem{proc}[theorem]{Procedure}

\numberwithin{equation}{section}

\newcommand{\bea}{\begin{eqnarray}}
\newcommand{\eea}{\end{eqnarray}}
\newcommand{\be}{\begin{eqnarray*}}
\newcommand{\ee}{\end{eqnarray*}}

\newcommand{\Z}{{\mathbb Z}}

\newcommand{\C}{{\mathbb C}}
\def\dsum#1#2{\mbox{$\sum\limits_{#1}^{#2}$}}
\def\II#1#2{{{[#1,#2]}}}
\newcommand{\fb}{{\mathfrak b}}
\newcommand{\fg}{{\mathfrak g}}
\newcommand{\fh}{{\mathfrak h}}

\newcommand{\Hom}{{\rm Hom}}

\def\HC#1{{\mathcal C}^{\rm Lexi}_{#1}}
\def\TC#1{{\mathcal C}^{\rm Trun}_{#1}}
\def\NC#1{{\mathcal C}^{\rm Norm}_{#1}}

\def\PED{}%{\hfill\mbox{$\Box$}}
\def\lm{_{\ssc \uparrow}}
\def\bm{^{\ssc \uparrow}}
\def\chiP#1{\chi^{\rm sum}(#1)}
\def\soe{\preccurlyeq}

\def\a{\alpha}

\def\ch{{\rm ch{\ssc\,}}}

\def\b{\beta}
\def\d{\delta}
\def\D{\Delta}
\def\g{\gamma}
\def\G{\Gamma}
\def\l{\lambda}
\def\L{\Lambda}
\def\si{\sigma}

\def\sc{\scriptstyle}
\def\ssc{\scriptscriptstyle}
\def\dis{\displaystyle}

\def\ol{\overline}
\def\ul{\underline}
\def\wt{\widetilde}
\def\wh{\widehat}

\def\sh#1#2#3{{#1}^\rho_{\cp\nn#3_{#2}}}
\def\shb#1#2#3{(#1)^\rho_{\cp\nn#3_{#2}}}

\def\Rar{\Longrightarrow}

\def\D{\Delta}

\def\Lra{\Longleftrightarrow}

\def\bs{\backslash}
\def\hs{\hspace*}
\def\rb{\raisebox}

\def\ni{\noindent}

\def\dd#1#2#3#4{{\widehat c}{\ssc\,}_{#4}(#2)}
\def\cc#1#2#3#4{c_{#4}(#2)}

\def\PP{P}
\def\Dr{P^{\rm Lex}_r}

\def\N{\mathbb{N}}
\def\Z{\mathbb{Z}}

\def\C{\mathbb{C}}

\def\gl{{\mathfrak {gl}}}

\def\es{\epsilon}
\def\deg#1{[#1]}

\def\cp#1{{\dot #1}}
\def\tt{{\rm aty}}
\def\ot{{\rm typ}}
\def\OT#1#2#3{\Set(\ot_{#2})}
\def\nn{{n}}
\def\mm{{m}}
\def\reg{^{\rm reg}}
\def\ITEM{\item[\mbox{$\circ$}]}
\def\RHO{\rho}
\def\fI{{{\textit{\textbf I{\,}}}}}
\def\pS{{\rm max}}
\def\tau{\pi}
\def\bi{{\ul{i}}}
\def\bj{{\ul{j}}}
\def\fS{S}
\def\fSa{\wt S}
\def\Sr{{\it Sym}}
\def\Set{{\rm Set}}

\def\CHIK{\chi^{\rm Kac}}
\def\equa#1#2{
\begin{equation}\label{#1}#2\end{equation}}
\def\equan#1#2{$$#2$$}
\numberwithin{equation}{section}
\begin{document}
%
%%%%%%%%%%%%%%%%%%%%%%%%%%%%%%%%%%%%%%%%%%%%%%%%%%%%%%%%
\title[Character and dimension formulae for ${\mathfrak{gl}}_{m|n}$]
{Character and dimension formulae for\\ 
general linear superalgebra}
\author[Yucai Su]{Yucai Su}
\address{School of Mathematics and Statistics,
University of Sydney, NSW 2006, Australia;\newline \indent
Department of Mathematics, Shanghai Jiaotong University, Shanghai
200030, China.} \email{yucai@maths.usyd.edu.au}
\author[R. B. ZHANG]{R. B. ZHANG}
\address{School of Mathematics and Statistics,
University of Sydney, NSW 2006, Australia}
\email{rzhang@maths.usyd.edu.au}

\begin{abstract}
The generalized Kazhdan-Lusztig polynomials for the finite
dimensional irreducible representations of the general linear
superalgebra are computed explicitly.
Using the result we
establish
a one to one correspondence
%an equivalence
between the
set of composition factors of an arbitrary $r$-fold atypical $\gl_{m|n}$-Kac-module
and  the
set of composition factors of some $r$-fold atypical $\gl_{r|r}$-Kac-module.
%
%category of $r$-fold atypical $\gl_{r|r}$-modules and
%the category of $\gl_{m|n}$-modules consisting of such objects that the
%composition factors all belong to an arbitrary $r$-fold atypical block.
%
The result of Kazhdan-Lusztig polynomials is also applied to
prove a conjectural character formula put forward by van der Jeugt
et al in the late 80s. We simplify this character formula to cast
it into the Kac-Weyl form, and derive from it a closed formula for
the dimension of any finite dimensional irreducible representation
of the general linear superalgebra.
\end{abstract}
\maketitle

%\tableofcontents
%\newpage

%%%%%%%%%%%%%%%%%%%%%%%%%%%%%%%%%%%%%%%%%%%%%%%%%%%%%%%%%%%%%%%%
%
\section{Introduction}
%
%%%%%%%%%%%%%%%%%%%%%%%%%%%%%%%%%%%%%%%%%%%%%%%%%%%%%%%%%%%%%%%%

Formal characters of finite dimensional irreducible
representations of complex simple Lie superalgebras encapsulate rich
information on the structure of the representations themselves. In
his foundational papers  \cite{Kac0, Kac, Kac1, Kac2} on Lie superalgebras, Kac
raised the problem of determining the formal characters of finite dimensional
irreducible representations of Lie superalgebras, and developed a
character formula for the so-called typical irreducible
representations.
However, the problem turned out to be quite hard
for the so-called atypical irreducible representations. In the
early 80s Bernstein and Leites \cite{BL} gave a formula for the general
linear superalgebra, which produces the correct formal characters
for the singly atypical irreducible representations \cite{VHKT0}, but fails for
the multiply atypical irreducibles (e.g., the trivial
representation is multiply atypical). Since then much further
research was done on the problem. For the orthosymplectic
superalgebra ${\mathfrak{osp}}_{2|2n}$, van der Jeugt \cite{V} constructed
a character formula for all finite dimensional irreducible
representations (which are necessarily singly atypical).  There
were also partial results and conjectures in other cases.
Particularly noteworthy is the conjectural character formula for
arbitrary finite dimensional irreducible representations of
${\mathfrak{gl}}_{m|n}$ put forward by
van der Jeugt, Hughes, King and Thierry-Mieg \cite{VHKT}, which was the result of
extensive research carried out by the authors over several years period.
Their formula withstood the tests of large scale computer
calculations for a wide range of irreducible representations.
However, the full problem of determining the formal characters of
finite dimensional irreducible representations of Lie
superalgebras remained open until 1995 when Serganova \cite{Se96, Se98} used a
combination of geometric and algebraic techniques to obtain a
general solution.

Serganova's approach was based on ideas from Kazhdan-Lusztig
theory. She introduced some generalized Kazhdan-Lusztig
polynomials, the values of which at $q=-1$ determine the formal
characters of finite dimensional irreducible representations of
Lie superalgebras. Serganova's work was further developed in the
papers \cite{Zou, VZ, B, B1, CZ}. Particularly important is the
work of Brundan, developed a very practicable algorithm for computing
the generalized Kazhdan-Lusztig polynomials, by using quantum group techniques.
This enables him to gain sufficient knowledge on the generalized Kazhdan-Lusztig
polynomials to prove  the conjecture of \cite{VZ} on the composition
factors of Kac-modules.

In this paper, we shall further investigate Brundan's algorithm
and implement it to compute the generalized Kazhdan-Lusztig
polynomials for the finite dimensional irreducible representations
of the general linear superalgebra. A closed formula is obtained
for the generalized Kazhdan-Lusztig polynomials, which is
essentially given in terms of the permutation group of the
atypical roots (see Theorem \ref{theo7.1} for details).

The formula for the generalized Kazhdan-Lusztig polynomials is
quite explicit and easy to apply. It leads to a
relatively explicit character formula
for all the finite
dimensional irreducible representations (see Theorem \ref{ours}).
By analysing this formula we prove that the conjecture of van der Jeugt et al
\cite{VHKT} holds true for all finite dimensional irreducible $\gl_{m|n}$-modules.

A general fact in the context of Lie
superalgebras is that the character formula constructed by using
Kazhdan-Lusztig theory is always in the form of an infinite sum.
This makes the character formula rather unwieldy to use for, e.g.,
determining dimensions of finite dimensional irreducible
representations. Therefore, it is highly desirable to sum up the infinite series
to cast the
 character
formula into the Kac-Weyl form. This is done in
Theorem \ref{theo7.2}.

Equipped with Theorem  \ref{theo7.2} we are able to work out the dimension of
any finite dimensional irreducible representation of
${\mathfrak{gl}}_{m|n}$, and the result is given by the closed
formula of Theorem  \ref{theo7.3}. In the special case of singly atypical
irreducible representations, our dimension formula reproduces what
one obtains from the Bernstein-Leites character formula \cite{V1}.

In proving Theorem \ref{theo7.1} on the Kazhdan-Lusztig polynomials we
have introduced the notion of heights of a weight with respect to its atypical roots.
This notion proves to be extremely useful. All the results in this paper
can be presented using this concept. In particular, the Kazhdan-Lusztig
polynomial $K_{\lambda, \mu}(q)$ depends only on the heights of $\lambda$
and $\mu$ with respect to their atypical roots. This latter fact
enables one to reduce the study of $r$-fold atypical
%indecomposable
$\gl_{m|n}$-Kac-modules
to the study of $r$-fold atypical
%indecomposable
$\gl_{r|r}$-Kac-modules.
We state this precisely in Theorem \ref{theo7.1'},
which establishes a one to one correspondence between the set of
composition factors of an arbitrary $r$-fold atypical $\gl_{m|n}$-Kac-module
and  the set of composition factors of some $r$-fold atypical $\gl_{r|r}$-Kac-module.
%
%category of $r$-fold atypical $\gl_{r|r}$-modules and
%the category of $\gl_{m|n}$-modules consisting of such objects that the
%composition factors all belong to an arbitrary $r$-fold atypical block.
%
%Let ${\mathcal E}^{m|n}_{\mathcal B}$
%be the full subcategory of $\gl_{m|n}$-modules consisting of such objects that the
%composition factors all belong to the $r$-fold atypical block ${\mathcal
%B}$. Theorem \ref{theo7.1'} establishes an equivalence between every
%such category and the category of $r$-fold atypical $\gl_{r|r}$-modules.

The organization of the paper is as follows. In Section 2 we present some
background material on ${\mathfrak {gl}}_{m|n}$, which will be used throughout
the paper. In Section 3 we investigate the generalized Kazhdan-Lusztig
polynomials for finite dimensional irreducible ${\mathfrak {gl}}_{m|n}$-modules.
This section contains two main results, Theorem \ref{theo7.1} and Theorem \ref{theo7.1'}.
While Theorem \ref{theo7.1} gives an explicit formula for the Kazhdan-Lusztig polynomials,
Theorem \ref{theo7.1'} establishes
%an equivalence between the
%category ${\mathcal E}^{m|n}_{\mathcal B}$
%of $\gl_{m|n}$-modules (${\mathcal B}$ is an $r$-fold atypical block)
%and the category of $r$-fold atypical $\gl_{r|r}$-modules.
a one to one correspondence between the set of
composition factors of an arbitrary $r$-fold atypical  Kac-module over $\gl_{m|n}$
and  the set of composition factors of some $r$-fold atypical Kac-module $\gl_{r|r}$.
In Section 4 we first use Theorem \ref{theo7.1} to prove the
conjectural character formula of van der Jeugt et al \cite{VHKT} (see Theorem
\ref{Conj-for}),
then we re-write the formula into the Kac-Weyl form (see Theorem \ref{theo7.2}).
Finally we derive from Theorem \ref{theo7.2} a closed formula for the
dimension of any finite dimensional irreducible representation of the
general linear superalgebra (see Theorem \ref{theo7.3}).

\section{Preliminaries}
\label{Preliminaries}
%%%%%%%%%%%%%%%%%%%%%%%%%%%%%%%%%%%%%%%%%%%%%%%%%%%%%%%%%%%%%%%%
%
%
\def\L{\lambda}\def\l{\mu}
We explain some basic notions of Lie superalgebras here and refer
to \cite{Kac, HKV, VHKT} for more details. We shall work over the
complex number field $\C$ throughout the paper. Given a
$\Z_2$-graded vector space $W=W_{\bar 0}\oplus W_{\bar 1}$, we
call $W_{\bar 0}$ and $W_{\bar 1}$ the even and odd subspaces,
respectively.
%The elements of $W_{\bar 0}\cup W_{\bar 1}$ will be called homogeneous.
Define a map $[\ ]: W_{\bar 0}\cup W_{\bar 1}
\rightarrow \Z_2$ by $[w] = {\bar i}$ if $w \in W_{\bar i}$. For
any two $\Z_2$-graded vector spaces $V$ and $W$, the space of
morphisms $\Hom_{\C}(V, W)$ (in the category of $\Z_2$-graded
vector spaces) is also $\Z_2$-graded with $\Hom_{\C}(V, W)_{\bar
k} =\sum_{\bar{i}+\bar{j}=\bar{k}}
\Hom_{\C}(V_{\bar i}, W_{\bar j})$. We write ${\rm End}_{\C}(V)$ for
$\Hom_{\C}(V, V)$.

Let $\C^{m|n}$ be the $\Z_2$-graded vector space with even subspace
$\C^m$ and odd subspace $\C^n$. Then ${\rm End}_{\C}(\C^{m|n})$ with
the $\Z_2$-graded commutator forms the {\it general linear superalgebra}.
To describe its structure, we choose a homogeneous basis
$\{ v_a \,|\, a\in \fI \}$,
for $\C^{m|n}$, where $\fI =\{1, 2, \ldots , m+n\}$, and $v_a$ is
even if $a\le m$, and odd otherwise. The general linear superalgebra relative to
this basis of $\C^{m|n}$ will be denoted by $\gl_{m|n}$,
which shall be further simplified to $\fg$ throughout the paper.
Let $E_{a b}$ be the matrix unit, namely, the $(m+n)\times(m+n)$-matrix
with all entries being zero except that at the $(a, b)$ position which is $1$.
Then $\{E_{a b}\, |\,a,b\in\fI \}$ forms a homogeneous basis of $\fg$,
with $E_{a b}$ being even if $a, b\le m$, or
$a, b> m$, and odd otherwise. For
convenience, we define the map $[\ ]: \fI \rightarrow \Z_2$ by
$[a]=\big\{\begin{array}{l l}
               \bar{0},  & \mbox{if} \ a\le m, \\
               \bar{1},  & \mbox{if} \ a>m.
              \end{array} $
Then the commutation relations of the Lie superalgebra can be written as
\be [E_{a b}, \ E_{c d}] &=& E_{a d}\delta_{b c}
 - (-1)^{([a]-[b])([c]-[d])} E_{c b}\delta_{a d}.
\ee

The upper triangular matrices form a {\it Borel subalgebra} $\fb$ of
$\fg$,
which contains the {\it Cartan subalgebra} $\fh$ of diagonal matrices.
Let $\{\es_a \,|\, a\in\fI \}$ be the basis
of $\fh^*$ such that
$\es_a(E_{b b})=\d_{a b}.$
The supertrace induces a {\em bilinear form} $(\;
, \: ): \fh^*\times \fh^* \rightarrow \C$ on $\fh^*$ such
that
$$(\es_a, \es_b)=(-1)^{[a]} \delta_{a b}.$$
Relative to the Borel subalgebra $\fb$,
the roots of $\fg$
can be expressed as $\es_a-\es_b, \:\, a\ne b$, where
$\es_a-\es_b$ is even if $[a]+[b] = \bar{0}$ and odd
otherwise. The set of the {\em positive roots}  is $\D^+=\{
\es_a-\es_b \,|\, a< b\}$, and the set of {\em simple roots}
is $\{\es_a-\es_{a+1} \,|\, a<m+n\}$.

We denote $\fI ^1=\{1,2,...,m\}$ and
$\fI ^2=\{1,2,...,n\}$. We also set $\d_\zeta=\es_{\cp\zeta}$
for $\zeta\in\fI ^2$, where
we use the notation
$$\cp \zeta=\zeta+m.
$$
%Then $\fI =\fI ^1\cup\fI ^2$.
Then the sets of {\em positive
even roots} and {\em odd roots} are respectively
\begin{eqnarray*}
\D_0^+&\!\!\!=\!\!\!&\{
\es_i-\es_j,\,
\d_\zeta-\d_\eta
\,|\,1\le i<j\le m,\, 1\le\zeta<\eta\le n\},
\\
\D_1^+&\!\!\!=\!\!\!&\{
\es_i-\d_\zeta\,|\,i\in\fI ^1,\,
\zeta\in\fI ^2\}.
\end{eqnarray*}
The Lie algebra $\fg$ admits a $\Z_2$-consistent $\Z$-grading
$$\fg=\fg_{-1}\oplus\fg_0\oplus\fg_{+1},
\mbox{ where }\fg_0=\fg_{\bar0}\cong\gl(m)\oplus\gl(n)
\mbox{  and }\fg_{\pm1}\subset\fg_{\bar1},
$$
with
$\fg_{+1}$ (resp.~$\fg_{-1}$) being the nilpotent subalgebra spanned by the odd positive (resp.~negative) root
spaces.
We define a total order on $\D_1^+$ by
%%%%%%%%%%%%%%%%%%%%%%%%%%%%%%%%%%%
\equa{0.0}
{
\es_i-\d_\zeta<\es_j-\d_\eta\,\ \ \ \Lra\,\ \ \
\zeta-i<\eta-j\mbox{ or }\zeta-i=\eta-j\mbox{ but }i>j.
}

An element in $\fh^*$ is called a {\em weight}.
A weight $\L\in\fh^*$ will be written in terms of the $\es\d$-basis as
%%%%%%%%%%%%%%%%%%%%%%%%%%%%%%%%%%%
\equa{weight1}
{\mbox{$
\L=(\L_1,...,\L_m\,|\,\L_{\cp 1},...,\L_{\cp n})=
\sum\limits_{i\in\fI ^1}\L_i\es_i-
\sum\limits_{\zeta\in\fI ^2}\L_{\cp\zeta}\d_\zeta,
$}}
where we have adopted an unusual (but convenient)
sign convention for $\L_i$'s.
Thus $\L_i=(\L,\es_i)$, called the {\it $i$-th entry of $\L$} for $i\in\fI ^1$, and
$\L_{\cp\zeta}=(\L,\d_\zeta)$, called the {\it $\cp\zeta$-th entry of $\L$}
for $\zeta\in\fI ^2.$
A weight $\L$ is
called
\begin{eqnarray}
\label{int}
\!\!\!\!\!\!\!\!
\mbox{\em integral}\ \ \,
&\Lra&
\L_i,\L_{\cp\zeta}\in\Z
\ \ \mbox{ for }\ i\in\fI ^1,\ \zeta\in\fI ^2;
\\
\label{dom}
\!\!\!\!\!\!\!\!
\mbox{\em dominant}
&\Lra&
{2(\Lambda, \, \alpha)}/{(\alpha,\, \alpha)}\ge0\mbox{ \ \ for all positive
even roots $\a$ of $\fg$, namely,}
\nonumber\\
\!\!\!\!\!\!\!\!
&&
\L_1\ge...\ge\L_m,\ \ \ \ \L_{\cp1}\le...\le\L_{\cp n}.
\end{eqnarray}
Denote by $\PP$ (resp.~$\PP_+$) the set of integral (resp.~dominant integral)
weights. Using notation (\ref{weight1}), $\PP$ coincides with the set of
the $m|n$-tuples of integers, thus $\PP$ is also denoted by $\Z^{m|n}$, and $\PP_+$ by
$\Z_+^{m|n}$.

Let $\rho_0$ (resp.~$\rho_1$) be half the sum of positive even
(resp.~odd) roots,
and let $\rho=\rho_0-\rho_1$. Then
%%%%%%%%%%%%%%%%%%%%%%%%%%%%%%%%%%%%%%%%%%%%%%%%%%%%%%%%
\begin{eqnarray}
%%%%%%%%%%%%%%%%%%%%%%%%%%%%%%%%%%%%%%%%%%%%%%%%%%%%%%%%
\label{rho0}
\!\!\!\!\!\!\!\!\rho_0
&\!\!\!\!=\!\!\!\!&\frac{1}{2}\big(\mbox{$\sum\limits_{i=1}^m$}(m-2i+1)\es_i+
\mbox{$\sum\limits_{\zeta=1}^n$}(n-2\zeta+1)\d_\zeta\big)
\nonumber\\[-4pt]
&\!\!\!=\!\!\!&\frac{1}{2}(m-1,m-3,...,1-m\,|\,1-n,3-n,...,n-1),
\\
\nonumber%\label{rho1}
\!\!\!\!\!\!\!\!\rho_1
&\!\!\!\!=\!\!\!\!&\frac{1}{2}\big(n\mbox{$\sum\limits_{i=1}^m$}\es_i
-m\mbox{$\sum\limits_{\zeta=1}^n$}\d_\zeta\big)
=
\frac{1}{2}(n,...,n\,|\,m,...,m),
\\
\nonumber%\label{rho}
\!\!\!\!\rho&\!\!\!\!=\!\!\!\!&
\rho'\!-\!\mbox{$\frac{m+n+1}{2}$}{\bf1},\mbox{ where }
\rho'\!=\!(m,...,2,1{\ssc\,}|{\ssc\,}1,2,...,n),\,{\bf1}=(1,...,1{\ssc\,}|{\ssc\,}1,...,1).
\end{eqnarray}
For all purposes, we can replace $\rho$ by $\rho'$. Therefore, from here on we shall denote
\equa{rho}
{
\RHO=(m,...,2,1\,|\,1,2,...,n).
}
Let $W=\Sr_m\times \Sr_n$ be the {\em Weyl group} of $\fg$,
where $\Sr_m$ is the {\it symmetric group of degree $m$}.
We define the {\em dot action} of $W$ on $\PP$
by
\equa{dot-action-of-W}
{w\cdot\l=w(\l+\rho)-\rho\mbox{ \ \ \ for \ \ }w\in W,\,\l\in \PP.
}
An integral weight $\L$ is called
\begin{enumerate}
\ITEM
 {\it regular}
or {\it non-vanishing} (in sense of \cite{HKV, VHKT}) if it is $W$-conjugate under the dot action to a
dominant weight
(which is denoted by $\L^+$ throughout the paper);
\ITEM
{\it vanishing\,} otherwise (since the right-hand side of
(\ref{typical-char}) is vanishing in this case, cf.~(\ref{vanishing})).
\end{enumerate}
Obviously,
\equa{l-ru}
{
\L\mbox{ \ is regular}\ \Lra\
\biggl\{
\begin{array}{ll}
\L_1+m,\,\L_2+m-1,\,...,\,\L_m+1&\mbox{are all distinct, and}
\\[4pt]
\L_{\cp1}+1,\,\L_{\cp2}+2,\,...,\,\L_{\cp n}+n&\mbox{are all distinct}.
\end{array}
}

Let $\L$  in (\ref{weight1}) be a regular weight.
%we define the {\it atypicality matrix of $\L$}
%to be the $m\times n$ matrix  $A(\L)=(A(\L)_{i,\zeta})_{m\times n}$, where
%(cf.~(\ref{rho}))
%%%%%%%%%%%%%%%%%%%%%%%%%%%%%%%%%%%
%\equan{Aty-matrix}
%{
%A(\L)_{i,\zeta}=(\L+\rho,\es_i-\d_\zeta)
%=(\L_i+m+1-i)-(\L_{\cp \zeta}+\zeta),\ i\in\fI ^1,\,\zeta\in\fI ^2.
%}
A positive odd root $\es_i-\d_\zeta$ is an {\it atypical root of $\L$} if
\equa{aty-root-L}
{(\L+\rho,\es_i-\d_\zeta)=(\L_i+m+1-i)-(\L_{\cp \zeta}+\zeta)=0.
}
Denote by $\G_\L$ the set of atypical roots of $\L$
(cf.~(\ref{(1)}) and (\ref{aty-roots})):
\equa{G-L}
{\G_\L=\{\es_i-\d_\zeta
\,|\,\,(\L+\rho,\es_i-\d_\zeta)=0\}.
}
Set $r=\#\G_\L.$
We also denote $\#\L=r$, called the {\it degree of
atypicality of $\L$}. A weight $\L$
is called
\begin{enumerate}
\ITEM
{\it typical} if $r=0$;
\ITEM
{\it atypical} if $r>0$ (in this case $\L$ is also called
an {\em $r$-fold atypical weight}).
\end{enumerate}

Let $V=\oplus_{\l\in\fh^*}V_\l$ be a {\em weight module} over $\fg$, where
\equan{weight-space}
{V_\l=\{v\in V\,|\,hv=\l(h)v,\,\forall\,h\in\fh\}\mbox{ \ with \ }\dim V_\l<\infty,
}
is the weight space of weight $\l$. The {\it character} $\ch V$ is defined to be
\equa{char-def}
{\mbox{$
\ch V=\sum\limits_{\l\in\fh^*}({\rm dim\,}V_\l) e^\l,
$}}
where $e^\l$ is the {\em formal exponential}, which will be regarded as an element
of an additive group isomorphic to $\fh^*$ under $\l\mapsto e^\l$.
Then $\ch V$ is an element of the {\em completed group algebra}
\equa{epsilon}
{\mbox{$
\varepsilon=\bigl\{\sum\limits_{\l\in\fh^*}a_\l e^\l\,\bigl|\,a_\l\in\C,\,
a_\l=0\mbox{ except $\l$ is in a finite union of ${\mathcal Q}_\L$}\bigr\},
$}}
where for $\L\in\fh^*$,
\equan{Q}
{\mbox{$
{\mathcal Q}_\L=\bigl\{\L-\sum\limits_{\a\in\D^+}
i_\a\a\in\fh^*\,\bigl|\,i_\a\in\Z_+\bigr\}.
$}}

For every integral dominant weight $\L$, we denote by $V^{(0)}(\L)$ the
finite-dimensional
irreducible $\fg_0$-module with highest weight $\L$. Extend it to a
$\fg_0\oplus\fg_{+1}$-module by putting $\fg_{+1}V^{(0)}(\L)=0$. Then
the {\em Kac-module} $\ol V(\L)$ is the induced
module
%%%%%%%%%%%%%%%%%%%%%%%%%%%%%%%%%%%
\equan{Kac-module}
{
\ol V(\L)={\rm Ind}_{\fg_0\oplus\fg_{+1}}^{\fg}V^{(0)}(\L)\cong U(\fg_{-1})
\otimes_{\C}V^{(0)}(\L).
}
Denote by $V(\L)$ the {\em irreducible
module} with highest weight $\L$ (which is the unique irreducible quotient
module of $\ol V(\L)$).

The following result is due to Kac \cite{Kac1, Kac2}:
%%%%%%%%%%%%%%%%%%%%%%%%%%%%%%%%%%%%%%%%%%%%%%%%%%%%%%%%
\begin{theorem}\label{typical-theo}
%%%%%%%%%%%%%%%%%%%%%%%%%%%%%%%%%%%%%%%%%%%%%%%%%%%%%%%%
If $\L$ is a dominant integral typical weight, then
$V(\L)=\ol V(\L)$, and
%%%%%%%%%%%%%%%%%%%%%%%%%%%%%%%%%%%
\equa{typical-char}
{
\ch V(\L)=\ch \ol V(\L)=\frac{L_1}{L_0}\,
\mbox{$\sum\limits_{w\in W}$}\,
\es(w)e^{w(\L+\rho)},
}
where
$\es(w)$ is the signature of $w\in W$, and
\equa{l-0}
{\mbox{$
L_0=\prod\limits_{\a\in\D_0^+}(e^{\a/2}-e^{-\a/2}),\,\ \
L_1=\prod\limits_{\b\in\D_1^+}(e^{\b/2}
+e^{-\b/2}).
$}}
\end{theorem}
%
%\begin{remark}\label{rema2.1}
%\rm
Since any finite-dimensional irreducible $\fg$-module
is either a typical module or
is a tensor module of $V(\L)$ with one-dimensional module
for some $\L\in \PP_+$,
in the rest of the paper
there is no loss of generality in restricting our attention to
integral weights $\L$.
%\end{remark}
%%%%%%%%%%%%%%%%%%%%%%%%%%%%%%%%%%%%%%%%%%%%%%%%%%%%%%%%%%%%%%%%
\section{Kazhdan-Lusztig polynomials}
\label{KL-polynomial}
%%%%%%%%%%%%%%%%%%%%%%%%%%%%%%%%%%%%%%%%%%%%%%%%%%%%%%%%%%%%%%%%
\def\ol{\bar}

\subsection{The height vectors and the $c$-relationship}
\label{sub3-c-related}
Let $\L\in\fh^*$ be an $r$-fold atypical dominant integral weight:
\equa{(1)}
{\mbox{$
\L\!=\!(\L_1,...,
\put(4,11){$\line(0,1){8}$}
\put(4,19){$\line(1,0){180}\,\rb{-4pt}{$\g_r$}\,\line(1,0){30}$}
%\put(4,-1){$\line(0,-1){5}$}
%\put(-53,-15){\footnotesize $r$-th atypical entry}
\L_{\mm_r}
, ...,
%\put(4,-1){$\line(0,-1){15}$}
%\put(-28,-25){\footnotesize typical entry}
\L_i, ...,
\put(4,10){$\line(0,1){5}$}
\put(4,15){$\line(1,0){20}\rb{-5pt}{$\,\g_1\,$}\line(1,0){30}$}
%\put(4,-1){$\line(0,-1){5}$}
%\put(-20,-15){\footnotesize $1$-st atypical entry}
\L_{\mm_1}
, ..., \L_m\,\bigl|\ \L_{\cp 1}, ...,
\put(4,10){$\line(0,1){5}$}
\put(4,15){$\line(-1,0){45}$}
%\put(4,-1){$\line(0,-1){5}$}
%\put(-40,-15){\footnotesize $1$-st atypical entry}
\L_{\cp \nn_1}
, ...,
%\put(4,-1){$\line(0,-1){15}$}
%\put(-28,-25){\footnotesize typical entry}
\L_{\cp\zeta}, ...,
\put(4,11){$\line(0,1){8}$}
\put(4,19){$\line(-1,0){40}$}
%\put(4,-1){$\line(0,-1){5}$}
%\put(-25,-15){\footnotesize $r$-th atypical entry}
\L_{\cp \nn_r}
,...,\L_{\cp n})$},
}
with the set of atypical roots
(cf.~(\ref{aty-root-L}), (\ref{G-L}) and (\ref{0.0}))
%%%%%%%%%%%%%%%%%%%%%%%%%%%%%%%%%%%
\equa{aty-roots}
{
\G_\L=\{\g_1,...,\g_r\}, \mbox{ \ \ \ where \ }
\g_1=\es_{\mm_1}-\d_{\nn_1}<...<\g_r=\es_{\mm_r}-\d_{\nn_r},
}
and $\mm_r<...<\mm_1,\ \nn_1<...<\nn_r$.
%Then we have
%\equa{aty-}
%{\g_1<...<\g_r
%\mbox{ \ \ \ (cf.~)}.
%}
We call $\g_s$ the {\it $s$-th atypical root of $\L$} for $s=1,...,r$.
For convenience, we introduce the notation $\L^\rho$ for the {\it $\rho$-translation
of $\L$}:
\equa{l-rho}
{
\L^\rho=\L+\rho.
}
Thus $\{(\mm_s,\nn_s)\,|\,s=1,...,r\}$ is the maximal set of pairs $(i,\zeta)$ satisfying
$\L^\rho_i=\L^\rho_{\cp\zeta}$ by (\ref{aty-root-L}).
We define the {\it atypical tuple of $\L$}
\equa{denote-t-f}
{\tt_\L=(\L^\rho_{\cp\nn_1},...,\L^\rho_{\cp\nn_r})=
(\L_{\cp\nn_1}+\nn_1,...,\L_{\cp\nn_r}+\nn_r)\in\Z^r,
}
and call the $s$-th entry of $\tt_\L$ the {\it $s$-th atypical entry of $\L$} for $s=1,...,r$.
We also define the {\it typical tuple of $\L$}
\equa{ot-f}
{
\ot_\L\in\Z^{m-r|n-r}
}
to be the element obtained
from $\L^\rho$ by deleting all entries $\L^\rho_{\mm_s},\,\L^\rho_{\cp\nn_s}$
for $s=1,...,r$. Thus all entries of $\ot_\L$,
called the {\it typical entries of $\L$}, are distinct by (\ref{l-ru}).
\begin{definition}\label{defi-height}
\rm
Corresponding to each atypical root $\g_s$ of $\L$, we
define the {\it $\g_s$-height of $\L$} (the {\it height of $\L$ with respect to the $s$-th atypical root})
\equa{s-height}
{
h_s(\L)=\L_{\mm_s}-\nn_s+s\mbox{ \ \ \ for \ \ }s=1,...,r.
}
We also introduce the {\it height vector of $\L$} and the {\it height of $\L$} respectively:
\equa{height}
{
h(\L)=(h_1(\L),...,h_r(\L)), \ \ \ \ \ |h(\L)|=\dsum{s=1}r h_s(\L).
}
\end{definition}
\begin{remark}\label{reme-height}
\rm
As we shall see later, the concept of heights of a weight with respect to its atypical
roots is extremely useful. In fact, the Kazhdan-Lusztig polynomials are completely determined
by the height vectors of the weights involved (see Theorem \ref{theo7.1} and
Theorem \ref{theo7.1'}).
\PED\end{remark}
\begin{example}\label{exam7.1-0}
\rm
Suppose $\L$ is the weight
\equa{exam1}
{
\L=\bigl(7,\stackrel{^{\sc4}}{{6}},5,\stackrel{^{\sc3}}{{5}},\stackrel{^{\sc2}}{3},3,2,
\stackrel{^{\sc1}}{2},0\,\bigl|\,\,1,\stackrel{^{\sc1}}{2},3,
\stackrel{^{\sc2}}{4},4,\stackrel{^{\sc3}}{{5}},\stackrel{^{\sc4}}{7},7\bigr)\in\Z^{9|8}_+,
}
where we put a label $t$ over a pair of entries to indicate
that they are the entries associated with the $t$-th
atypical root.
Then
\begin{eqnarray}
\label{exam1'}
\L^\rho&\!\!\!=\!\!\!&\bigl(16,\ul{14},12,\ul{11},
\ul{8},7,5,
\ul{4},1\,\bigl|\,\,2,\ul{4},6,
\ul{8},9,\ul{11},\ul{14},15\bigr),
\\
\label{exam1'''-}
\tt_\L&\!\!\!=\!\!\!&(4,8,11,14),
\\
\label{exam1''}
\ot_\L&\!\!\!=\!\!\!&(16,12,7,5,1\,|\,2,6,9,15),
\\\label{exam-height}
h(\L)&\!\!\!=\!\!\!&(1,1,2,3),
\end{eqnarray}
where the underlined integers in $\L^\rho$ are entries of atypical tuple $\tt_\L$.
\end{example}
We define the partial order `$\soe$' on $P_+$ for $\L,\l\in P_+$ by
\equa{g-prec-f}
{
\l\soe\L\ \ \ \Lra\ \ \ \#\l=\#\L,\ \,\tt_\l\le \tt_\L\mbox{ \ and \ }
\ot_\l=\ot_\L,
}
where the partial order ``$\le$'' on $\Z^r$
is defined for $a=(a_1,...,a_r),\,b=(b_1,...,b_r)
\in\Z^r$ by
\equa{partial-order-on-Z-r}
{
a\le b\ \ \ \Lra\ \ \
a_i\le b_i,\ \ \forall\,i\in\II{1}{r}.
}

Now suppose $\l$ is another $r$-fold atypical dominant integral
weight with atypical roots:
\equa{mm'-s}
{
\g'_1=\es_{\mm'_1}-\d_{\nn'_1}< ...< \g'_r=\es_{\mm'_r}-\d_{\nn'_r}.
}
We define (cf.~(\ref{prop-1-1}) below)
\equa{l-g-f'}
{
\ell_s(\L,\l)=h_s(\L)-h_s(\l)=
\L_{\mm_s}-\l_{\mm'_s}+\nn'_s-\nn_s \ \ \ \mbox{ for \ }s\in\II{1}{r},
}
and define
the {\it length between $\L$ and $\l$} to be
\equa{l-g-f}
{
\ell(\L,\l)=\dsum{s=1}r\ell_s(\L,\l)=|h(\L)|-|h(\l)|.
}
(In general $\ell(\L,\l)$ is not necessarily non-negative, but when $\l\soe\L$,
it is indeed non-negative.)
%Let $\l\in P_r$ be another weight with atypical roots $\g'_s=
%\es_{\mm'_s}-\d_{\nn'_s}$, $1\le s\le r$, being as in (\ref{mm'-s}).
\begin{remark}\label{rem?}
The height $|h(\L)|$ of $\L$ turns out to be
the absolute length $\ell(\L)$ defined by
Brundan in \cite[\S3-g]{B}, and the length $\ell(\L,\l)$
coincides with the length $\ell(\l,\L)$ defined in \cite[\S3-g]{B}.
\end{remark}
For $1\le s\le t\le r$,
we denote
\equa{dis-s-t}
{
d_{s,t}(\L)=h_t(\L)-h_s(\L)
=\L_{m_t}-\L_{m_s}-\nn_t+\nn_s+t-s.
}
The fact that $d_{s,t}(\L)$ is non-negative is not obvious, but one can observe that
it is the number of integers between the $s$-th atypical entry
$\L^\rho_{\cp\nn_s}$ and the $t$-th atypical entry
$\L^\rho_{\cp\nn_t}$ which are not entries of $\L^\rho$, namely (cf.~(\ref{set-R(f)}), (\ref{I-st})
and (\ref{prop-1-1}))
\equa{d(s,t)}
{d_{s,t}(\L)=
\#\bigl(\II{\L^\rho_{\cp\nn_s}}{\L^\rho_{\cp\nn_t}}\bs \Set(\L^\rho)\bigr).
}
%which follows from the more general result in Proposition \ref{prop-1},
Because of this fact, we call $d_{s,t}(\L)$
the {\it distance between two atypical roots $\g_s$ and $\g_t$ of $\L$}.
Here and below we use the notation
\equa{set-R(f)}
{\Set(\l)=\mbox{ the set of the entries of a weight }\l,
}
and the notation
\equa{I-st}
{
\II{i}{j}=\biggl\{\begin{array}{cl}
\{k\in\Z\,|\,i\le k\le j\}&\mbox{if }i\le j,
\\[4pt]
\emptyset&\mbox{otherwise},\end{array}
\mbox{ \ \ \ for \ \ }i,j\in\Z.
}
(There will be no danger of confusing this notation with the Lie bracket as the
later will not be used in the remainder of the paper.)
%We remark that this notation will not be confused
%with the Lie bracket since we do not need to
%use that in the rest of the paper.

One can generalize (\ref{d(s,t)}) to obtain the following proposition,
which will be used in the proof of Theorem \ref{theo7.1'}.
\begin{proposition}\label{prop-1}
Let $\L,\l$ be $r$-fold atypical weights with $\l\soe\L$. For any $s,t\in\II{1}{r}$
%$($not necessarily $s\le t)$
with
$\l^\rho_{\cp\nn'_s}\le \L^\rho_{\cp\nn_t}$, we have
\equa{prop-1-1}
{
\#\bigl(\II{\l^\rho_{\cp\nn'_s}+1}{\L^\rho_{\cp\nn_t}}\bs \Set(\ot_\L)\bigr)
=h_t(\L)-h_s(\l)+t-s.
}
\end{proposition}
\begin{proof}
By (\ref{s-height}) and the fact that $\L^\rho_{\mm_s}=\L^\rho_{\cp\nn_s}$,
the right-hand side of (\ref{prop-1-1}) is equal to
\equan{Fact1}
{
\L_{\cp\nn_t}+\mm_t-(\l_{\cp\nn'_s}+\mm'_s)+2(t-s)
=\L^\rho_{\cp\nn_t}-\l^\rho_{\cp\nn'_s}-(\mm'_s-\mm_t+\nn_t-\nn'_s)-2(s-t).
}
Thus (\ref{prop-1-1}) is equivalent to
\equa{Fact2}
{
\#\bigl(\II{\l^\rho_{\cp\nn'_s}+1}{\L^\rho_{\cp\nn_t}}\cap \Set(\ot_\L)\bigr)
=\mm'_s-\mm_t+\nn_t-\nn'_s+2(s-t).
}
Note from (\ref{g-prec-f}) that $\ot_\L=\ot_\l$.
By (\ref{dom}), (\ref{l-ru}) and (\ref{(1)}), in $\L^\rho$ (resp., $\l^\rho$),
the number of typical entries to the left of entry $\L^\rho_{\cp\mm_t}$
(resp., $\l^\rho_{\mm'_s}$) is $\mm_t-1-r+t$ (resp., $\mm'_s-1-r+s$).
Thus the number of typical entries
to the left of the entry $\l^\rho_{\mm'_s}$ which are
in $\II{\l^\rho_{\mm'_s}+1}{\L^\rho_{\mm_t}}$
is $\mm'_s-\mm_t+s-t$.
Similarly, the number of typical entries
to the right of the entry $\l^\rho_{\mm'_s}$ which are
in $\II{\l^\rho_{\cp\nn'_s}+1}{\L^\rho_{\cp\nn_t}}$
is $\nn_t-\nn'_s+s-t$. Hence we obtain (\ref{Fact2}).
\end{proof}
For any two distinct
atypical roots $\gamma_s$ and $\gamma_t$, there exist
unique positive even roots $\a_{st}$ and $\b_{st}$ such that
$(\a_{st}, \b_{st})=0$, and the composite of the actions of the
Weyl group elements $\sigma_{\a_{st}}$, $\sigma_{\b_{st}}$, respectively
corresponding to reflections with respect to these even roots, send $\gamma_s$
to $\gamma_t$, namely, $\gamma_t =\sigma_{\a_{st}}
\sigma_{\b_{st}}(\gamma_s)$.
In terms of explicit formulae, we have
\equa{a-st}
{
\a_{st}=\es_{\mm_t}-\es_{\mm_s},\ \ \ \b_{st}=\d_{\nn_s}-\d_{\nn_t}\in\D_0^+
\mbox{ \ \ \ for \ \ }1\le s<t\le r.
}
The concept of $c$-relationship defined below was first introduced in \cite{HKV}
from a different point of view.
\begin{definition}\label{c-related}
\rm
For $s\le t$, two atypical roots $\g_s,\g_t$ of $\L$ are called
{\it $c$-related} (in the sense of \cite{HKV})
or {\it connected} (in the sense of \cite{VZ})
if $s=t$ or $d_{s,t}(\L)<t-s$. The later relation has the nice interpretation that
\equan{Inter}{
\begin{array}{l}
\mbox{\it the distance $($i.e., $d_{s,t}(\L))$ between two atypical roots is smaller than}
\\
\mbox{\it the number $($i.e., $t-s)$ of atypical roots between them.}
\end{array}
}
Relation $d_{s,t}(\L)<t-s$ is also equivalent to
\equa{(3)}
{\L_{m_t}-\L_{m_s}<\nn_t-\nn_s,
}
or in terms of weights,
\equan{(3)'}{
\frac{2(\L,\a_{st})}{(\a_{st},\a_{st})}
+\frac{2(\L,\b_{st})}{(\b_{st},\b_{st})}
<
\frac{2(\rho,\a_{st})}{(\a_{st},\a_{st})}
+\frac{2(\rho,\b_{st})}{(\b_{st},\b_{st})}.
\eqno(\ref{(3)})'}
We define
\begin{eqnarray}
\label{c'-st}
&&\!\!\!\!\!\!\!\!\!\!\!\!\!\!\!\!\!\!\!\!\!\!\!\!
\cc^\L_{s,t}=\biggl\{\begin{array}{ll}
1&\mbox{if the
atypical roots $\g_s,\g_t$ of $\L$ are $c$-related},\\
0&\mbox{otherwise}.
\end{array}
\end{eqnarray}
We also define
\begin{eqnarray}
\label{define-d-st}
&&\!\!\!\!\!\!\!\!\!\!\!\!
\dd^\L_{s,t}=\biggl\{\begin{array}{ll}
1&\mbox{if \ }\cc^\L_{s,s}=\cc^\L_{s,s+1}=...=\cc^\L_{s,t}=1,\\
0&\mbox{otherwise},
\end{array}
\end{eqnarray}
and we say that {\it $\g_s,\g_t$ are strongly $c$-related} if
$\dd^\L_{s,t}=1$. \PED
%namely, all atypical roots $\g_p$'s of $\L$ with
%$s\le p\le t$ are $c$-related to the atypical root $\g_s$ of $\L$.
%Two atypical roots $\g_s,\g_t$ of $\L$ are called
%{\it disconnected} if $\cc^\L_{s,t}=0$.
\end{definition}

Note that $c$-relationship is reflexive
and transitive but not symmetric
($\cc^\L_{t,s}$ is not defined when $t>s$).

\begin{example}\label{exam7.1}
\rm
Let $\L$ be the weight in (\ref{exam1}).
By (\ref{(3)}),
\begin{enumerate}
\item[] $\cc^\L_{1,2}=1$
since $\L_{\mm_2}-\L_{\mm_1}=3-2<4-3=\nn_2-\nn_1$;
\item[]
$\cc^\L_{1,3}=1$ since $\L_{\mm_3}-\L_{\mm_1}=5-2<6-2=\nn_3-\nn_1$;
\item[]
$\cc^\L_{1,4}=1$ since $\L_{\mm_4}-\L_{\mm_1}=6-2<7-2=\nn_4-\nn_1$;
and
\item[]
$\cc^\L_{s,t}=0$ for any other pair $(s,t)$.
\end{enumerate}
Thus
\equan{exam1.1}
{\mbox{
$\dd^\L_{1,2}=\dd^\L_{1,3}=\dd^\L_{1,4}=1$ and
$\dd^\L_{s,t}=0$ for any other pair $(s,t).$
}
%\eqno\PED
}
\end{example}
\begin{remark}\label{reme-in-D-r}
\rm
If $\L$ is regular but not necessarily dominant,
we shall generalize the notions $c,\wh c,d,h$ to $\L$
by defining them with respect to the dominant weight $\L^+$. For instance, $\cc^\L_{s,t}=\cc^{\L^+}_{s,t}$.
Sometimes even if $\L$ is not regular but lexical $($in the sense of Definition
$\ref{defi-lexical})$,
one can still define the $c$-relationship by using the %language of
distance (\ref{dis-s-t}).
\PED
\end{remark}

%The following discussion may illustrate the significance of the concept of
%$c$-relation\-ship.

\subsection{Lexical weights}
Let $\L$
%given in (\ref{(1)})
be an $r$-fold atypical regular weight
(not necessarily dominant)
with the set
$\G_\L=\{\g_1,...,\g_r\}$
of atypical roots
ordered according to (\ref{aty-roots}).
%(cf.~(\ref{l-ru})).
We call $\L$ {\it lexical} if its atypical tuple $\tt_\L$ is lexical in the following sense:
\begin{definition}\label{defi-lexical}
\rm
An element $a=(a_1,...,a_r)\in\Z^r$ is called {\it lexical}
if
\equa{lexical}
{
a_1\le ...\le a_r.
}
%\PED
\end{definition}

The two sets $P_r$ and $\Dr$ to be defined below will be frequently used throughout this
section.
\begin{definition}\label{defi-two-sets}
\rm
We denote by $P_r$ the set of $r$-fold atypical regular weights $\L$
of the form  (\ref{(1)}) such that the atypical roots of $\L$
can be ordered as in (\ref{aty-roots})
% satisfying (\ref{aty-})
and that the typical tuple $\ot_\L\in\Z_+^{m-r|n-r}$ is
dominant as a weight for $\gl_{m-r|n-r}$.

We denote by $\Dr$ the subset of $P_r$ consisting of the lexical weights of $P_r$.
\PED
\end{definition}
%\begin{convention}
%\label{conv}
%\rm
%We usually use the superscript `$\L$' to indicate that
%a notation is associated with $\L$, like in $\cc^\L_{s,t},$ $\dd^\L_{s,t}$.
%However, when confusion is unlikely
%to occur, the superscript will be dropped, like in $k_s$ defined below.
%\end{convention}

\subsection{The $r$-tuple %$(k_1,...,k_r)$
of positive integers associated with $\L$}
The $r$-tuple associated with $\L\in \Dr$ defined below
 was first introduced in \cite{VZ}.
%seems to be quite useful.

\begin{definition}\label{k-i-L}
\rm
Define the $r$-tuple $(k_1,...,k_r)$ of positive integers associated with $\L\in \Dr$
in the following way:
%for each $s=r,...,1$,
each $k_s$ is the smallest positive integer such that
\equa{cond7.1}
{
(\L+\theta_t k_t\g_t)+k_s\g_s\mbox{ \ \
is regular for all \ }t=s+1,s+2,...,r\mbox{ and }\theta_t\in\{0,1\}.
%\vs{-6pt}
}
%\ \vs{-6pt}\PED
\end{definition}
The following lemma gives a way %tells how
to compute $k_s$'s.
%To state the lemma, we need to introduce one more notation.
%
For $\L\in\Dr$ and $s\in\II{1}{r}$, we set
\equa{p-s}
{
\pS^\L_s=\max\{p\in
%\{1, 2, \cdots, r\}
\II{s}{r}
\,|\,
%p\ge s, \,
\dd^\L_{s,p}=1\}
}
to be the maximal number $p\in\II{s}{r}$ satisfying
the condition that the $p$-th
atypical root $\g_p$ of $\L$ is strongly $c$-related to $\g_s$.
One immediately sees that
\equa{c-s-c-t}
{
\pS^\L_t\le \pS^\L_s\mbox{ \ \ \ for any $t$ with }s\le t\le \pS^\L_s.
}
\begin{lemma}\label{lemm7.1}
Let $\L\in \Dr$.
\begin{enumerate}
\item
\label{lemm7.1-1}
For $s\in\II{1}{r}$, $k_s$ is the integer such that
$\L^\rho_{\cp\nn_s}+k$ is the \mbox{$(\pS^\L_s+1-s)$-th} smallest integer bigger than $\L^\rho_{\cp\nn_s}$ and
not in the entry set $\Set(\L^\rho)$, i.e.,
\equa{k-i}
{k_s=\min\bigl\{k>0\,\bigl|\ \#\bigl(\II{\L^\rho_{\cp\nn_s}}{\L^\rho_{\cp\nn_s}+k}
\bs \Set(\L^\rho)\bigr)
=\pS^\L_s+1-s\bigr\}.
}
\item
\label{lemm7.1-2}
The tuple $(k_r,...,k_1)$ is the lexicographically
smallest tuple of positive integers such that for all $\theta=(\theta_1,...,
\theta_r)\in\{0,1\}^r$, $\L+\sum_{s=1}^r\theta_sk_s\g_s$ is regular. Thus
$(k_r,...,k_1)$ is the tuple satisfying \cite[{\rm Main Theorem}]{B}.
\end{enumerate}
\end{lemma}
\begin{proof} (\ref{lemm7.1-1})
Denote by $k'_s$ the right-hand side of (\ref{k-i}).
Obviously, $k'_r$ is the smallest positive
integer such that $\L+k'_r\g_r$ is
regular because for any $0<k<k'_r$,
by definition the integer $\L^\rho_{\cp\nn_r}+k$ which is equal
to $\L^\rho_{\mm_r}+k$ already
appears in the entry set $\Set(\L^\rho)$ and thus
$\L+k'_r\g_r$ is not regular by (\ref{l-ru}).

If $s<r$, by induction on $\pS^\L_s$, it is straightforward to see
that
$\L^\rho_{\cp\nn_s}+k'_s$ is the smallest integer (bigger than $\L^\rho_{\cp\nn_s}$)
which is not in
$\Set(\L^\rho+\theta_tk'_t\g_t)$
for $\theta_t\in\{0,1\}$ and $s< t\le r$. Thus
$k'_s$ is the smallest positive integer satisfying (\ref{cond7.1}).

(\ref{lemm7.1-2})
Similarly,
$\L^\rho_{\cp\nn_s}+k'_s$ is also the smallest integer
(bigger than $\L^\rho_{\cp\nn_s}$) which is not in
$\Set(\L^\rho+\sum_{t=s}^r\theta_tk'_t\g_t)$
for $\theta=(\theta_s,...,
\theta_r)\in\{0,1\}^{r+1-s}$, i.e,
$(k'_r,...,k'_s)$ is the lexicographically
smallest tuple of positive integers such that for all $\theta=(\theta_s,...,
\theta_r)\in\{0,1\}^{r+1-s}$, the weight
$\L+\sum_{t=s}^r\theta_tk'_t\g_t$ is regular.
\end{proof}
Lemma \ref{lemm7.1}(\ref{lemm7.1-1}) allows us to compute $k_s$ by the following
procedure.
\begin{proc}\label{proc1}
First set $S=\Set(\L^\rho)$. Suppose we have computed $k_r,...,k_{s+1}$.
To compute $k_s$, we count the numbers in the set $S$ starting with $\L^\rho_{\cp\nn_s}$
until we find a number, say $k$, not in $S$. Then $k_s=k-\L^\rho_{\cp\nn_s}$.
Now add $k$ into the set $S$, and continue.
\end{proc}
\begin{example}\label{exam-k-i}
\rm
Let $\L$ be given in $(\ref{exam1})$.
Using the above procedure we obtain
$(k_4,k_3,k_2,k_1)=(3,2,2,14)$ (cf.~(\ref{exam1'})).
\PED\end{example}
\begin{remark}\label{reme-k-i}
\rm
If $\L\in P_r$ (not necessarily in $\Dr$),
we can still compute $k_s$ by the above procedure, but the difference lies in that the
$k_s$'s are computed not in the order $s=r,...,1$, but in the order that
each time we compute $k_s$ with $\L^\rho_{\cp\nn_s}$ being the largest among all those
$\L^\rho_{\cp\nn_s}$'s, the corresponding $k_s$'s of which are not yet computed.
\PED\end{remark}
\subsection{Raising operators}
Following \cite{B}, we define the {\it raising operator}
 $R_{\mm_s,\cp\nn_s}$ on $P_r$ by
\equa{raising}
{
R_{\mm_s,\cp\nn_s}(\L)=\L+k_s\g_s\ \ \ \mbox{ for \ }\L\in P_r
\mbox{ \ and \ }s\in\II{1}{r},
}
where $\g_s,k_s$ are defined in (\ref{aty-roots}) and Definition \ref{k-i-L} respectively.
Obviously, for $\L,\l\in \Dr$,
\equa{fg-ol-r}
{\l=R_{\mm_s,\cp\nn_s}(\L)\ \ \Rar\ \ \ot_\L=\ot_\l
\mbox{ \ \ (cf.~(\ref{ot-f}))}.
}
Denote $\N=\{0,1,...\}$, %the set of {\it natural numbers},
and let
$\theta=(\theta_1,...,\theta_r)\in\N^r$.
We define
\equa{R'-theta-f}
{R'_\theta(\L)=(R^{\theta_1}_{\mm_1,\cp\nn_1}\circ\cdots\circ R^{\theta_r}_{\mm_r,\cp\nn_r}(\L))^+,
}
where in general $\l^+$ denotes the unique dominant element
which is $W$-conjugate under the dot action to $\l$ (cf.~(\ref{dot-action-of-W})).

Let $\l\in P_r$ be another weight with atypical roots $\g'_s=
\es_{\mm'_s}-\d_{\nn'_s}$, $1\le s\le r$, being as in (\ref{mm'-s}).
Then we have (cf.~\cite[\S3-f]{B})
\equa{prec}
{\l\soe\L\ \ \ \Lra\ \ \
\#\l=\#\L=:r,\mbox{ \ \ and \ }\exists\,\theta\in\N^r\mbox{ \ with \ }
R'_\theta(\l)=\L.
}
For convenience, we denote
\equa{g--f}
{
\l\prec\!\prec\L\mbox{ \ \ \ if \ }\ot_\l=\ot_\L\mbox{ \ and \ }
\max\{\l^\rho_{\cp\nn'_s}\,|\,s\in\II{1}{r}\}\le
\min\{\L^\rho_{\cp\nn_s}\,|\,s\in\II{1}{r}\}.
}
\subsection{Definitions of $S^\L$ and $\fS^{\L,\l}$}
The symmetric group $\Sr_r$ of degree $r$ acts on
$\Z^r$ by permuting entries. This action
induces an action on $P_r$
given by
\equa{(4)}{{\sc\!\!\!\!\!\!\!\!\!}
\si(\L)\!=\!(\L_1,{\sc...\,},
\put(4,-1){$\line(0,-1){5}$}
\put(-25,-13){\small atypical entries permuted}
\L_{\mm_{\si(r)}},{\sc...\,},\L_i,{\sc...\,},
\put(4,-1){$\line(0,-1){5}$}
\L_{\mm_{\si(1)}},{\sc...\,},\L_m\,|\,
\L_{\cp 1},{\sc...\,},
\put(4,-1){$\line(0,-1){5}$}
\put(-25,-13){\small atypical entries permuted}
\L_{\cp \nn_{\si(1)}},{\sc...\,},\L_{\cp\zeta},{\sc...\,},
\put(4,-1){$\line(0,-1){5}$}
\L_{\cp \nn_{\si(r)}},{\sc...\,},\L_{\cp n}),\!\!\!\!\!\!\!
}
for $\si\in \Sr_r$ and $\L\in P_r$.
With this action on $P_r$, the group $\Sr_r$ can be regarded as a subgroup of
$W$, such that every element is of even parity.
Thus we also have the {\it dot action}
\equa{Action-of-s-r}
{\si\cdot\L=\si(\L+\rho)-\rho
\mbox{ \ \ \ for \ }\si\in \Sr_r.
}
\begin{notation}\label{Notation2}
\rm
Let $\L,\l\in \Dr$. Define
$S^\L$ to be the subset of the symmetric group $\Sr_r$ consisting of permutations
$\si$ which do not change the order of $s<t$ when the atypical roots
$\g_s$ and $\g_t$ of $\L$ are strongly $c$-related. That is,
\equa{s-f-r}{
S^\L=\{\si\in \Sr_r\,|\,\,
\si^{-1}(s)<\si^{-1}(t)\mbox{ \ for \ all \ }s<t
\mbox{ \ with \ }\dd^\L_{s,t}=1\},
}
where $\dd^\L_{s,t}$ is defined in (\ref{define-d-st}).
We also define $\fS^{\L,\l}$ to be subset of $S^\L$ consisting of permutations $\si$
such that $\l\soe\si\cdot\L$, namely
\equa{s-f-g}{
\fS^{\L,\l}=\{\si\in S^\L\,|\,\,\l\soe\si\cdot\L\}.
%\vs{-6pt}
}
%\vs{-6pt}\PED
\end{notation}
For convenience, we denote
\equa{s-f-g'}{
\fSa^{\L,\l}=\{\si\in \Sr_r\,|\,\,\l\soe\si\cdot\L\}.
}
Thus $\fS^{\L,\l}=S^\L\cap \fSa^{\L,\l}.$
\begin{example}\label{exam-s-f}
\rm
If $\L$ is the weight in $(\ref{exam1})$, then
$S^\L=\{\si\in \Sr_4\,|\,\si(1)=1\}\cong \Sr_3, $
is a subgroup of $\Sr_4$ $($however in general $S^\L$ is not a subgroup$)$.
\PED\end{example}
Let $\ell(\si)$ denote the normal {\it length function on $\Sr_r$}, namely
\begin{eqnarray}
\label{ell-si}
&\!\!\!\!\!\!&\!\!\!\!\!\!
\ell(\si)=\dsum{s=1}r\ell(\si,s), \mbox{ \ \ where}
\\
\label{ell(si,s)}
&\!\!\!\!\!\!&\!\!\!\!\!\!
\ell(\si,s)=\#\{t>s\,|\,\si(t)<\si(s)\}\mbox{ \ \ \ for \ \ }s=1,...,r.
\end{eqnarray}
For any subset $B\subset \Sr_r$, we define the {\it $q$-length function
of $B$} by:
\equa{l-q-function}
{
B(q)=\mbox{$\sum\limits_{\si\in B}$}q^{\ell(\si)}.
}
\begin{proposition}\label{lemm7.1+}
Let $\L,\l\in \Dr$ with $\l\soe\L$.
We have
\begin{eqnarray}
\label{Sr(q)}
%\nonumber
\Sr_r(q)&\!\!\!\!=\!\!\!\!&(1+q)(1+q+q^2)\cdots(1+q+\cdots+q^{r-1})=
\mbox{$\prod\limits_{s=1}^r$}\frac{q^s-1}{q-1},
\\[-8pt]
\label{s-f(q)}
%\nonumber
S^\L(q)&\!\!\!\!=\!\!\!\!&
\mbox{$\prod\limits_{s=1}^r$}(q^s-1)/
\mbox{$\prod\limits_{s=1}^r$}(q^{\pS^\L_s-s+1}-1),
\\
\label{s-fg-(q)}
%\nonumber
\fSa^{\L,\l}(q)&\!\!\!\!=\!\!\!\!&
(1+q+\cdots+q^{r-i_r})(1+q+\cdots+q^{r-1-i_{r-1}})\cdots
\nonumber
\\[-4pt]
&\!\!\!\!=\!\!\!\!&
\mbox{$\prod\limits_{s=1}^r$}\frac{q^{s+1-i_s}-1}{q-1},
\end{eqnarray}
where
\equa{i-s}{
i_s=\min\{i\in\II{1}{r}\,|\,
\l^\rho_{\cp\nn'_s}\le\L^\rho_{\cp\nn_i}\}
\mbox{ \ \ \  for \ }s\in\II{1}{r}.
}
\end{proposition}
\begin{proof}
First we compute $\fSa^{\L,\l}(q)$.
Elements $\si\in \fSa^{\L,\l}$ can be easily described as follows
(cf.~Example \ref{exam7.4}):
for each $s=r,...,1$,
suppose for all $t>s$, $\si(t)$ have been chosen, then
$\si(s)$ can be any of $s-i_s+1$ integers $r,r-1,...,i_s$ which have not yet been
taken by the $\si(t)$'s for $t>s$. We order the elements of $\{r,...,i_s\}\bs\{\si(r),...,\si(s+1)\}$
in descending order, and denote by
\equa{order-1}{
\{r,...,i_s\}\bs\{\si(r),...,\si(s+1)\}=\{x_1>...>x_{s-i_s+1}\}.
}
Then each choice of $\si(s)=x_k$ for $k=1,...,s-i_s+1$
contributes $k-1$ to the length $\ell(\si)$.
Thus we have (\ref{s-fg-(q)}). Since $\Sr_r=\fSa^{\L,\l}$ for any
$\L,\l\in \Dr$ with $\l\prec\!\prec\L$ (cf.~(\ref{s-f-g'})), and
when $\l\prec\!\prec\L$, all $i_s$'s are equal to $1$, we obtain
(\ref{Sr(q)}), which is a well-known formula.

Consider $S^\L$ defined in (\ref{s-f-r}), which can be re-written as
\equan{s-f-r'}
{
S^\L
=\{\si\in \Sr_r\,|\,
\si^{-1}(s)<\si^{-1}(t)\mbox{ for all }s,t \mbox{ with }s<t\le \pS^\L_s\}.
}
Since for each $s=1,...,r$, we cannot change the order of $s$ and $t$ for
$s<t\le \pS^\L_s$, we shall remove the factor $1+q+\cdots+q^{\pS^\L_s-s}$
from $\Sr_r(q)$. Thus we obtain (\ref{s-f(q)}).
\end{proof}
Similar to (\ref{s-fg-(q)}), we also have
\\[8pt]\hs{0ex}
$\dis
\fSa^{\L,\l}(q)\!=\!
\mbox{$\prod\limits_{s=1}^r$}\frac{q^{s+1-j_s}\!-\!1}{q-1},
\mbox{ where }
j_s=\max\{j\in\II{1}{r}\,|\,
\L^\rho_{\cp\nn_s}\ge \l^\rho_{\cp\nn'_j}\}\mbox{  for }s\in\II{1}{r}.
$\hfill(\ref{s-fg-(q)})$'$
\subsection{The $q$-length function of $\fS^{\L,\l}$}
Elements $\si\in \fS^{\L,\l}$ can be described in the following way
(cf.~the proof of Proposition \ref{lemm7.1+} and Example \ref{exam7.4}):
\begin{desc}\label{desc-1}
\rm
For $s=r,...,1$,
each $\si(s)$ can be any one of the numbers $r,r-1,...,i_s$ which has not yet been
occupied by $\si(t)$ for some $t>s$ (cf.~(\ref{order-1})), with an additional condition that
if $\dd^\L_{a,b}=1$ for some $a<b$ such that $b$ has not yet been chosen, then
$\si(s)\ne a$.
\PED\end{desc}
We can associate each $\si\in \fS^{\L,\l}$ with
a graph defined as follows:
Put $r$ weighted points at the bottom such that the $s$-th point
(which will be referred
to as {\it point $s^-$}) has weight $
%w(s^-)=
\l^\rho_{\cp\nn'_s}$.
Similarly we put $r$ weighted points on the top such that the $s$-th point
(which will be referred to as {\it point $s^+$}) has weight
$
%w(s^+)=
\L^\rho_{\cp\nn_s}$.
Two points $s^+$ and $t^+$ on the top are connected by a line if and only if
$\dd^\L_{s,t}=1$, in this case we say that the two points are {\it linked.}
Note that if $s^+$ is linked to $t^+$, then
$s^+$ is linked to $p^+$ for all $p$ with $s<p<t$ (cf.~(\ref{c-s-c-t})).
Each $\si\in \Sr_r$ can be represented by a graph which is obtained by drawing a line
(denoted by $L(s)$) between the bottom point $s^-$ and the top
point $\si(s)^+$ for each $s\in\II{1}{r}$.
The length $\ell(\si)$ is simply the number of crossings %\linebreak[4]
\mbox{\rb{9pt}{$\put(0,0){\line(3,-2){20}}\put(0,0){\phantom{\line(1,0){20}}}
\put(20,0){\line(-3,-2){20}}$} \hs{4ex}} in the graph.
Then $\si\in \fS^{\L,\l}$ if and only if
the weight of $\si(s)^+$ is not less than the weight of $s^-$ for each $s$, and
in the case $L(t)$ crosses $L(s)$, the two points $\si(s)^+$ and $\si(t)^+$
on the top cannot be linked, i.e., a graph with the part \rb{9pt}{$\put(0,0){\line(3,-2){20}}\put(0,0){\line(1,0){20}}
\put(20,0){\line(-3,-2){20}}$} \hs{4ex} is not allowed.
See the example below.
\begin{example}\label{exam7.4}
\rm
Let $\L$ be given in $(\ref{exam1})$ and
\begin{eqnarray}
\label{exam1.3-1}
&&\!\!\!\!\!\!\!\!\!\!\!\!
\l\,=\,\bigl(\,\,7,\,\,\,4,\,\,\,\stackrel{^{\sc4}}{{4}},\,\,\,\stackrel{^{\sc3}}{{4}},{\sc\,}2,{\sc\,}1,
\stackrel{^{\sc2}}{1},
\stackrel{^{\sc1}}{1},0\,
\bigl|\,
1,\,\stackrel{^{\sc1}}{1},\stackrel{^{\sc2}}{1},2,
4,\,\,\stackrel{^{\sc3}}{{4}},\,\,\stackrel{^{\sc4}}{4},\,\,7\bigr),
\mbox{ \ \ and so}\\
\label{exam1.3-1'}
&&\!\!\!\!\!\!\!\!\!\!\!\!
\l^\rho=\bigl(16,12,\ul{11},\ul{10},7,5,
\ul{4},\ul{3},1
\,\bigl|\,\,
2,\ul{3},\ul{4},6,
9,\ul{10},\ul{11},15\bigr)
\mbox{ \ \ \ (cf.~(\ref{exam1'}))}.
\end{eqnarray}
Thus $\l\soe\L$. The
elements of $\fSa^{\L,\l}$ correspond to the following graphs:
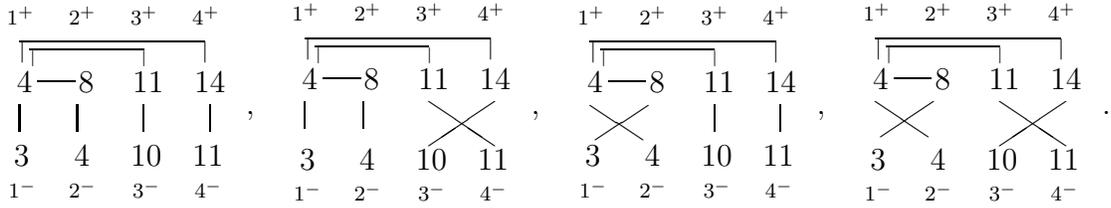
\equan{exam7.4-1}
{
\begin{array}{llll}
\put(15,4){\line(1,0){14}}
\put(12,16){\put(0,-5){$\ssc|$}\line(1,0){43}\put(-1,-5){$\ssc|$}}
\put(8,19){\put(0,-6){$\sc|$}\line(1,0){70}\put(-1,-6){$\sc|$}}
\put(8,-5){\line(0,-1){10}}\put(30,-5){\line(0,-1){10}}\put(55,-5){\line(0,-1){10}}
\put(80,-5){\line(0,-1){10}}
\stackrel{^{^{^{^{^{^{\sc1^+}}}}}}}{\dis4}& \stackrel{^{^{^{^{^{^{\sc2^+}}}}}}}{\dis8}&
\stackrel{^{^{^{^{^{^{\sc3^+}}}}}}}{\dis11}&\stackrel{^{^{^{^{^{^{\sc4^+}}}}}}}{\dis14}
\\[15pt]
\,\,\stackrel{^{{\dis3}}}{\sc1^-}& \stackrel{^{{\dis4}}}{\sc2^-}&
\stackrel{^{{\dis10}}}{\sc3^-}& \stackrel{^{{\dis11}}}{\sc4^-}
\end{array}
\,,
\
\begin{array}{llll}
\put(15,4){\line(1,0){14}}
\put(12,16){\put(0,-5){$\ssc|$}\line(1,0){43}\put(-1,-5){$\ssc|$}}
\put(8,19){\put(0,-6){$\sc|$}\line(1,0){70}\put(-1,-6){$\sc|$}}
\put(8,-5){\line(0,-1){10}}\put(30,-5){\line(0,-1){10}}
\put(55,-5){\line(3,-2){25}}
\put(80,-5){\line(-3,-2){25}}
\stackrel{^{^{^{^{^{^{\sc1^+}}}}}}}{\dis4}& \stackrel{^{^{^{^{^{^{\sc2^+}}}}}}}{\dis8}&
\stackrel{^{^{^{^{^{^{\sc3^+}}}}}}}{\dis11}&\stackrel{^{^{^{^{^{^{\sc4^+}}}}}}}{\dis14}
\\[15pt]
\,\,\stackrel{^{{\dis3}}}{\sc1^-}& \stackrel{^{{\dis4}}}{\sc2^-}&
\stackrel{^{{\dis10}}}{\sc3^-}& \stackrel{^{{\dis11}}}{\sc4^-}
\end{array}
\,, \
\begin{array}{llll}
\put(15,4){\line(1,0){14}}
\put(12,16){\put(0,-5){$\ssc|$}\line(1,0){43}\put(-1,-5){$\ssc|$}}
\put(8,19){\put(0,-6){$\sc|$}\line(1,0){70}\put(-1,-6){$\sc|$}}
\put(8,-5){\line(3,-2){20}}\put(30,-5){\line(-3,-2){20}}
\put(55,-5){\line(0,-1){10}}
\put(80,-5){\line(0,-1){10}}
\stackrel{^{^{^{^{^{^{\sc1^+}}}}}}}{\dis4}& \stackrel{^{^{^{^{^{^{\sc2^+}}}}}}}{\dis8}&
\stackrel{^{^{^{^{^{^{\sc3^+}}}}}}}{\dis11}&\stackrel{^{^{^{^{^{^{\sc4^+}}}}}}}{\dis14}
\\[15pt]
\,\,\stackrel{^{{\dis3}}}{\sc1^-}& \stackrel{^{{\dis4}}}{\sc2^-}&
\stackrel{^{{\dis10}}}{\sc3^-}& \stackrel{^{{\dis11}}}{\sc4^-}
\end{array}
\,, \
\begin{array}{llll}
\put(15,4){\line(1,0){14}}
\put(12,16){\put(0,-5){$\ssc|$}\line(1,0){43}\put(-1,-5){$\ssc|$}}
\put(8,19){\put(0,-6){$\sc|$}\line(1,0){70}\put(-1,-6){$\sc|$}}
\put(8,-5){\line(3,-2){20}}\put(30,-5){\line(-3,-2){20}}
\put(55,-5){\line(3,-2){25}}
\put(80,-5){\line(-3,-2){25}}
\stackrel{^{^{^{^{^{^{\sc1^+}}}}}}}{\dis4}& \stackrel{^{^{^{^{^{^{\sc2^+}}}}}}}{\dis8}&
\stackrel{^{^{^{^{^{^{\sc3^+}}}}}}}{\dis11}&\stackrel{^{^{^{^{^{^{\sc4^+}}}}}}}{\dis14}
\\[15pt]
\,\,\stackrel{^{{\dis3}}}{\sc1^-}& \stackrel{^{{\dis4}}}{\sc2^-}&
\stackrel{^{{\dis10}}}{\sc3^-}& \stackrel{^{{\dis11}}}{\sc4^-}
\end{array}\, .\
}
We see that $i_4=i_3=3,\,i_2=i_1=1$
$($where $i_s$'s are defined in $(\ref{i-s}))$, thus $\fSa^{\L,\l}(q)=(1+q^{4-3})1(1+q^{2-1})1=(1+q)^2$,
which agrees with the above graphs.

The elements of $\fS^{\L,\l}$ are represented by the first two of the above graphs. Thus
$$
\fS^{\L,\l}(q)=1+q.
%\eqno\PED
$$
\end{example}
\def\ZF{Z}
Now let us compute the $q$-length function $\fS^{\L,\l}(q)$.
First we introduce a family of $q$-functions $\ZF_q(x;b)$ defined on the set of
pairs $(x,b)$ of lexical $r$-tuples $x=(x_1,...,x_r)$, $b=(b_1,...,b_r)\in\Z^r$
satisfying $1\le b_s\le s$ for $s\in\II{1}{r}$, i.e.,
\equa{z-r-elements}
{
x_1\le x_2\le ...\le x_r,\ \ \
b_1\le b_2\le\cdots\le b_r
\mbox{ \ and \ }1\le b_s\le s, \ \forall\,s\in\II{1}{r}.
}
\begin{definition}\label{defi-F}
\rm
Define the $q$-function $\ZF_q(x;b)$ as follows:
Set $\ZF_q(x;b)=0$ if (\ref{z-r-elements}) is not satisfied, and define
$\ZF_q(x;b)$ inductively on $r$ by:
\begin{eqnarray}
\label{F-par}
&&\!\!\!\!\!\!\!\!\!\!\!\!\!\!\!\!\!\!\!
\ZF_q(x_1;1)=1,\ \ \ \
\ZF_q(x_1,x_2;1,b_2)=1+\theta(x_2-x_1-1)\theta(1-b_2)q.
\\
\label{Define-F}
&&\!\!\!\!\!\!\!\!\!\!\!\!\!\!\!\!\!\!\!
\ZF_q(x;b)=
\ZF_q(x_1,...,x_{r-1};b^{(r-1)})
\nonumber\\
&&\!\!\!\!\!\!\!\!\!\!\!\!\!\!\!\!\!\!\!
\phantom{\ZF_q(x;b)=}
+\dsum{i=b_r}{r-1}
\theta(x_{i+1}{\sc\!}-{\sc\!}x_i{\sc\!}-{\sc\!}1)
\ZF_q(x_1,...,x_{i-1},x_{i+1}{\sc\!}-{\sc\!}1,...,x_r{\sc\!}-{\sc\!}1;b^{(r-1)})q^{r-i},
\end{eqnarray}
where $b^{(r-1)}=(b_1,...,b_{r-1})$, and $\theta(x)$ is the {\it step function} defined
by
\equa{d-function}
{
\theta(x)=\biggl\{\begin{array}{ll}
1&\mbox{if \ }x\ge0,\\
0&\mbox{otherwise}.
\end{array}
%\vs{-10pt}
}
%\vs{-4pt}\PED
\end{definition}
Note that there are only finite many functions $\ZF_q(x;b)$ for each fixed $r$. To see this,
for any lexical $x\in\Z^r$ and for
$1\le s\le t\le r$,
we define $\cc^x_{s,t}$ and $\dd^x_{s,t}$ analogous to (\ref{c'-st}) and (\ref{define-d-st})
by
\begin{eqnarray}
\label{c'-st-x}
&&\!\!\!\!\!\!\!\!\!\!\!\!\!\!\!\!
\cc^x_{s,t}=\biggl\{\begin{array}{ll}
1&\mbox{if \ }t=s\mbox{ or }x_t-x_s<t-s,\\
0&\mbox{otherwise},
\end{array}
\\
\label{define-d-st-x}
&&\!\!\!\!\!\!\!\!\!\!\!\!\!\!\!\!
\dd^x_{s,t}=\biggl\{\begin{array}{ll}
1&\mbox{if \ }\cc^x_{s,s}=\cc^x_{s,s+1}=...=\cc^x_{s,t}=1,\\
0&\mbox{otherwise},
\end{array}
\end{eqnarray}
and define $\wh x=(\wh x_1,...,\wh x_r)$ with
\equa{x'}
{
\wh x_s=\#\{p\in\II{1}{s-1}\,\,|\,\,\dd^x_{p,s}=0\}\mbox{ \ \ \ for \ } s\in\II{1}{r}.
}
Then one can prove that $\wh x$ is lexical and $0\le \wh x_s<s$ for $s\in\II{1}{r}$, and
\equa{F-x'}
{
\ZF_q(x;b)=\ZF_q(\wh x;b).
}
Thus
there are only finite many functions $\ZF_q(x,b)$ for each fixed $r$.

\begin{proposition}\label{q-length-S}
Let $\L,\l\in\Dr$ with $\l\soe\L$.
Denote
\equa{where}
{%\biggl\{
\begin{array}{l}
b^{\L,\l}=(i_1,...,i_r)
\in\Z^r,\mbox{  where  }
i_s=\min\{i\in\II{1}{r}\,|\,\L^\rho_{\cp\nn_s}\le\l^\rho_{\cp\nn'_i}\}
\end{array}
 \mbox{ for  $s\in\II{1}{r}$},
}
namely
$i_s$ is defined in $(\ref{i-s})$.
Let $h(\L)$ be the height vector of $\L$ defined in $(\ref{height}).$
Then the $q$-length function of $\fS^{\L,\l}$ is given by
\equa{q-s-Ll-function}
{
\fS^{\L,\l}(q)=\ZF_q(h(\L);b^{\L,\l})
%=\ZF_q(\wh x^\L;b^{\L,\l})
.
}
\end{proposition}
\begin{proof}
The proposition follows from
Description \ref{desc-1} and Definition \ref{c-related}.
\end{proof}
\subsection{Generalized Kazhdan-Lusztig polynomials}
Let $K_{\L,\l}(q)$ be the Kazhdan-Lusztig polynomial defined in \cite{Se98}
(which was denoted as $K_{\l,\L}(q)$ in \cite{B}, but we prefer its original notation in \cite{Se98}),
and
$\ell_{\L,\l}(q)$ the polynomial defined in \cite{B}
(which was denoted as $\ell_{\l,\L}(q)$ in \cite{B}).
Then \cite[Corollary 3.39]{B} states
that
\equa{ell-gf}{\mbox{$
K_{\L,\l}(q)=
\ell_{\L,\l}(-q^{-1})=\sum\limits_{\theta\in\N^{\#\L}:\:R'_\theta(\l)=\L}q^{|\theta|},
\mbox{ \ \ where }|\theta|=\dsum{s=1}r\theta_s,
$}}
for dominant weights $\L,\l$.
Let $\ell(\L,\l)$ be the length defined in (\ref{l-g-f}).
First let us see an example of computing $\ell(\L,\l)$.
\begin{example}\label{lenth-of-f-g}
\rm
Suppose $\L$ is given in $(\ref{exam1})$ and $\l$ given in $(\ref{exam1.3-1})$.
By (\ref{l-g-f'}),
\equan{exam1.3-3}
{
\begin{array}{lll}
\ell_1(\L,\l)=2-1+2-2=1,&&
\ell_2(\L,\l)=3-1+3-4=1,
\\[4pt]
\ell_3(\L,\l)=5-4+6-6=1,&&
\ell_4(\L,\l)=6-4+7-7=2.
\end{array}
}
Thus
$\ell(\L,\l)=7.$
\PED\end{example}

One of the main results of this paper is the following.
%%%%%%%%%%%%%%%%%%%%%%%%%%%%%%%%%%%%%%%%%%%%%%%%%%%%%%%%
\begin{theorem}\label{theo7.1}
Suppose $\L,\l$ are dominant weights.
Then $K_{\L,\l}(q)\ne0$ if and only if $\l\soe\L$ and in this case
\equa{k-l-p-f}{\mbox{$
K_{\L,\l}(q)
=
q^{\ell(\L,\l)}
\sum\limits_{\si\in \fS^{\L,\l}}q^{-2\ell(\si)}
=q^{\ell(\L,\l)}\fS^{\L,\l}(q^{-2}),
$}}
where $\fS^{\L,\l}$ is defined by $(\ref{s-f-g})$, and
$\fS^{\L,\l}(q)$ is the $q$-length function of
$\fS^{\L,\l}$ defined in $(\ref{l-q-function})$
and can be determined by $(\ref{q-s-Ll-function})$.
\end{theorem}
A weight $\L$ is said to be
\begin{enumerate}
\ITEM
{\it totally disconnected\,}
if $\cc^\L_{s,t}=0$ for all pairs $(s,t)$ with $s<t$;
\ITEM
{\it totally connected\,}
if $\cc^\L_{s,t}=1$ for all pairs $(s,t)$ with $s\le t$.
\end{enumerate}
From Theorem \ref{theo7.1}, one immediately obtains
\begin{coro}
\label{coro7.1}
Let $\L,\l\in \Dr$ with $\l\soe\L$, and let $S^\L(q),\,\fSa^{\L,\l}(q)$ be as in
Proposition $\ref{lemm7.1+}$. We have
\begin{enumerate}
\item %{\rm(1)}
If $\L$ is totally connected, then $K_{\L,\l}(q)=q^{\ell(\L,\l)}$.
\item %\par{\rm(2)}
If $\L$ is totally disconnected, then
$
K_{\L,\l}(q)=q^{\ell(\L,\l)}\fSa^{\L,\l}(q^{-2}).
$
\item %\par{\rm(3)}
If $\l\prec\!\prec\L$ $($recall~$(\ref{g--f}))$, then
$K_{\L,\l}(q)=q^{\ell(\L,\l)}S^\L(q^{-2}).$
\end{enumerate}
\end{coro}

\subsection{Proof of Theorem \ref{theo7.1}}
%We shall break the proof of Theorem \ref{theo7.1} into three lemmas.
%
%For convenience, we should always write an element $a$ in $\Z^r$ as $a=(a_1,...,a_r)$.
%In particular, the atypical tuple of $\nu$ (cf.~(\ref{aty-root-L})) is written as
%\equan{s-h}
%{
%\tt_\nu=\bigl(\sh{\nu}{1}{},...,\sh{\nu}{r}{}\bigr)\mbox{ \ \ \ for \ }\nu\in P_r.
%}
Using notations as above, we set (cf.~(\ref{ell-gf}))
%Set
$$\Theta^\L_\l=\{\theta\in\N^r\,|\,R'_\theta(\l)=\L\},
\mbox{ \ \ \ where \ }r=\#\L.$$
In the three lemmas below, we
shall establish a bijection %one to one correspondence
between $\Theta^\L_\l$ and $\fS^{\L,\l}$.
%This is achieved in  the following lemmas.
Theorem \ref{theo7.1} is then a simple consequence of this bijection.

First let us define a map
%\begin{lemma}\label{t-l-1}
%There exists a map
\equa{injection}{
\Theta^\L_\l\to \fSa^{\L,\l}:\ \theta\mapsto\si_\theta
}
%\end{lemma}
%\begin{proof}
in the following way.
Suppose $\theta\in\Theta^\L_\l$, i.e. (cf.~(\ref{mm'-s}) for notations $\mm'_s,\nn'_s$),
\equa{R'-theta}
{
\L=R'_\theta(\l)=(R^{\theta_1}_{\mm'_1,\cp\nn'_1}\circ\cdots\circ R^{\theta_r}_{\mm'_r,\cp\nn'_r}(\l))^+.
}
We denote
\equa{denote-gi}
{
\l^{(r)}=\l,\ \
\l^{(s-1)}=R^{\theta_s}_{\mm'_s,\cp\nn'_s}(\l^{(s)}) \
\mbox{ \ \ for \ } s=r,r-1,...,1.
}
Then each $\l^{(s-1)}$ is obtained from $\l^{(s)}$ by adding some number,
say, $K_s$, to its $\mm'_s$-th, $\cp\nn'_s$-th entries. In fact (we use the notation
$k_s^{\l^{(s,i)}}$ to denote the integer $k_s$ defined in (\ref{cond7.1})
with $\L$ replaced by $\l^{(s,i)}$)
\equa{K-s}
{\mbox{$
K_s=\sum\limits_{i=1}^{\theta_s}k_s^{\l^{(s,i)}},
\mbox{ \ \
where }\l^{(s,0)}=\l^{(s)},\ \l^{(s,i)}=R_{\mm'_s,\cp\nn'_s}(\l^{(s,i-1)}).
$}}
Thus
\equa{l-(0)}
{\l^{(0)}=\l+\dsum{s=1}r K_s\g'_s.
}
Since $\L=(\l^{(0)})^+$ and both $\L$ and $\l^{(0)}$ are regular,
there exists
a $\si\in \Sr_r$ uniquely
determined by $\theta$, such that (cf.~(\ref{(4)}) and (\ref{Action-of-s-r}))
\equa{g1-si-f}
{
\tt_{\l^{(0)}}=\si(\tt_\L)=\tt_{\si\cdot\L},\mbox{ \ \ i.e., \ }
%\l^{(0)}_{\cp\nn_s}+\nn_s=\L_{\cp\nn_{\si(s)}}+\nn_{\si(s)}
\shb{\l^{(0)}}{s}{}=\sh{\L}{\si(s)}{}
\mbox{ \ for  }s\in\II{1}{r}
\mbox{ \ (cf.~(\ref{denote-t-f}))}.
}
By (\ref{l-(0)}), $\tt_\l\le \tt_{\l^{(0)}}=\tt_{\si\cdot\L}$,
and by (\ref{g-prec-f}) and (\ref{prec}), $\ot_\l=\ot_\L$.
But $\ot_\L=\ot_{\si\cdot\L}$ (typical tuples are invariant under the dot action of
$\Sr_r$, cf.~(\ref{ot-f}), (\ref{(4)}) and (\ref{Action-of-s-r})).
Thus $\l\soe\si\cdot\L$ by (\ref{g-prec-f}),
i.e., $\si\in \fSa^{\L,\l}$ by definition (\ref{s-f-g'}).
We denote $\si$ by $\si_\theta$. Thus we obtain the map (\ref{injection}).
%\end{proof}
\begin{lemma}\label{t-l-2}
The map $(\ref{injection})$ is an injection. More precisely, suppose
$\si\in \Sr_r$ such that $\si=\si_\theta$ for some $\theta=(\theta_1,...,\theta_r)\in\Theta^\L_\l$,
then such $\theta$ is unique and is given by
\begin{eqnarray}
\label{theta-s}
\!\!\!\!&\!\!\!\!&\!\!\!\!
\theta_s=\theta'_s-2\ell(\si,s)
\mbox{ \ \ \ $($cf.~$(\ref{ell(si,s)}))$,\ \ where}
\\
\label{Q-s}
\!\!\!\!&\!\!\!\!&\!\!\!\!
\theta'_s=\#Q_s,\ \ Q_s=
\II{\sh{\l}{s}{'}+1}{\sh{\L}{\si(s)}{}}
\bs \OT (\L)\mbox{ \ \ \ for \ }s=r,r-1,...,1.
\end{eqnarray}
Thus by $(\ref{prop-1-1}))$, $\theta'_s$ is in fact $\theta'_s=h_{\si(s)}(\L)-h_s(\l)+\si(s)-s$.
\end{lemma}
\begin{proof}
The proof of the lemma is divided into the following cases.
\vskip4pt
{\it Case 1: $s=r$}.

We want to prove
$\theta_r=\theta'_r$.
Note that each time when we apply $R_{\mm'_r,\cp\nn'_r}$ to $\l$, the
$r$-th entry
$\sh{\l}{r}{'}$
 of the atypical tuple $\tt_\l$
reaches an integer in the set $Q_r$ (cf.~(\ref{k-i}) and (\ref{raising})), and
no integer in this set $Q_r$ can be skipped. Thus after applying $\theta'_r$ times,
this entry reaches the integer
%$\L_{\cp\nn_{\si(r)}}+\nn_{\si(r)}$
$\sh{\L}{\si(r)}{}$
(it is an entry of $\tt_\L$, thus not in the typical entry set $\OT(\L)$), and $\l$ becomes $\l^{(r-1)}$ (cf.~(\ref{denote-gi})). Thus $\theta_r=\theta'_r$.

One can also use the following arguments to prove
$\theta_r=\theta'_r$:
The integer
%$\L_{\cp\nn_{\si(r)}}+\nn_{\si(r)}$
$\sh{\L}{\si(r)}{}$
is the $r$-th entry of
the atypical tuple $\tt_{\l^{(r-1)}}$, which by definitions
(\ref{raising}) and (\ref{R'-theta-f}) and by (\ref{k-i}), is equal to the
$\theta_r$-th smallest integer bigger than
%$\l_{\cp\nn'_r}+\nn'_r$
$\sh{\l}{r}{'}$
and not in $\OT(\L)$. But by the definition of $\theta'_r$ in (\ref{Q-s}),
%$\L_{\cp\nn_{\si(s)}}+\nn_{\si(s)}$
$\sh{\L}{\si(s)}{}$
 is the $\theta'_r$-th
smallest integer bigger than
%$\l_{\cp\nn'_r}+\nn'_r$
$\sh{\l}{r}{'}$
and not in $\OT(\L)$. Thus $\theta_r=\theta'_r$.
\vskip4pt
{\it Case 2: $s=r-1$ and $\si(r-1)<\si(r)$}.

We want to prove $\theta_{r-1}=\theta'_{r-1}$.
There are two possibilities to consider.
\vskip4pt
{\it Subcase 2.i: Suppose
%$\l_{\cp\nn'_{r-1}}+\nn'_{r-1}=\L_{\cp\nn_{\si(r-1)}}+\nn_{\si(r-1)}$}.
$\sh{\l}{r-1}{'}=\sh{\L}{\si(r-1)}{}$}.
Then $\l^{(r-2)}=\l^{(r-1)}$ and
we obviously have $\theta_{r-1}=0=\theta'_{r-1}$. We are done.
\vskip4pt
{\it Subcase 2.ii: Suppose
%$\l_{\cp\nn'_{r-1}}+\nn'_{r-1}<\L_{\cp\nn_{\si(r-1)}}+\nn_{\si(r-1)}$}.
$\sh{\l}{r-1}{'}<\sh{\L}{\si(r-1)}{}$}.
Note that
$
%\L_{\cp\nn_{\si(r-1)}}+\nn_{\si(r-1)}
%<\L_{\cp\nn_{\si(r)}}+\nn_{\si(r)}.
\sh{\L}{\si(r-1)}{}
<\sh{\L}{\si(r)}{}.
$
Also note that
\equa{fact1-}
{\Set((\l^{(r-1)})^\rho)=(\Set(\l^\rho)\bs\{\sh{\l}{r}{'}\})
\cup\{\sh{\L}{\si(r)}{}\},
}
(recall the $\rho$-translated notation in (\ref{l-rho})
and note that the only difference between $\l^{(r-1)}$
and $\l$ is their $r$-th atypical entries).
We observe that
$\sh{\L}{\si(r-1)}{}$ is not
an entry of the typical tuple $\ot_\l=\ot_\L$ because of the regularity of $\L$, it is not
the $t$-th entry of the atypical tuple $\tt_\l$ either for any $t\le r-1$
because of the assumption that $\sh{\l}{r-1}{'}<\sh{\L}{\si(r-1)}{}$.
Thus $\sh{\L}{\si(r-1)}{}$ is not in (\ref{fact1-}).
But
it is in the set
\equa{fact2-}
{\II{\shb{\l^{(r-1)}}{r-1}{'}}{\shb{\l^{(r-1)}}{r}{'}}
=\II{\sh{\l
%^{(r-1)}
}{r-1}{'}}{\sh{\L}{\si(r)}{}}.
}
Hence there is at least an integer in (\ref{fact2-}) but not in (\ref{fact1-}). By
Definition \ref{c-related},
the $(r-1)$-th and $r$-th atypical roots
of $\l^{(r-1)}$
are not $c$-related.

Similarly, the $(r-1)$-th and $r$-th atypical roots
of $\l^{(r-1,i)}$ (cf.~(\ref{K-s}))
are not $c$-related
for all $i$ with $1\le i<\theta'_{r-1}$.
Thus by the arguments in Case 1, we need to apply the raising operator
$R_{\mm'_{r-1},\cp\nn'_{r-1}}$ to $\l^{(r-1)}$
exactly $\theta'_{r-1}$ times
in order to obtain $\l^{(r-2)}$. So $\theta'_{r-1}=\theta_{r-1}$.
\vskip4pt
{\it Case 3: $s=r-1$ and $\si(r-1)>\si(r)$}.

We want to prove $\theta_{r-1}=\theta'_{r-1}-2$.
Note that
\begin{eqnarray}
\label{entry-r-1}
\!\!\!\!
\shb{\l^{(r-1)}}{r-1}{'}
=%&\!\!\!=\!\!\!&
\sh{\l}{r-1}{'}
<\sh{\l}{r}{'}
\le \shb{\l^{(0)}}{r}{'}
=%\nonumber\\&\!\!\!=\!\!\!&
\sh{\L}{\si(r)}{}
<\sh{\L}{\si(r-1)}{}\ \mbox{ (cf.~(\ref{g1-si-f}))},
\end{eqnarray}
(recall that $\L,\l$ are dominant).
Thus we need to apply $R_{\mm'_{r-1},\cp\nn'_{r-1}}$ to $\l^{(r-1)}$ at least once.

Suppose after applying $i$ times of
$R_{\mm'_{r-1},\cp\nn'_{r-1}}$ to $\l^{(r-1)}$ ($i$ can be zero),
the \mbox{$(r-1)$-th} entry $\shb{\l^{(r-1)}}{r-1}{'}$
of the atypical tuple $\tt_{\l^{(r-1)}}$
reaches an integer, say $p$,
such that
\equa{such-that-}
{p<\sh{\L}{\si(r)}{}
\mbox{ \ \ but \ \ }
\II{p}{\sh{\L}{\si(r)}{}}\subset \Set(\L^\rho).
}
By (\ref{entry-r-1}), such $i$ must exist
since the $(r-1)$-th entry will finally reach the integer
$\sh{\L}{\si(r-1)}{}$.

Note that $p$ and $\sh{\L}{\si(r)}{}$ are respectively
the $(r-1)$-th and $r$-th entries of $\tt_{\l^{(r-1,i)}}$ (cf.~(\ref{K-s})). Thus by
(\ref{such-that-}) and Definition \ref{c-related},
the $(r-1)$-th and $r$-th atypical roots
of $\l^{(r-1,i)}$ are $c$-related.
Then by (\ref{k-i}) and (\ref{raising}),
when we apply $R_{\mm'_{r-1},\cp\nn'_{r-1}}$ to $\l^{(r-1,i)}$, the
$(r-1)$-th entry of $\tt_{\l^{(r-1,i)}}$
reaches an integer in the set
\equa{the-set}
{
\II{\sh{\L}{\si(r)}{}+1}{\sh{\L}{\si(r-1)}{}}
\bs\OT ({\l^{(r-1,i)}})
=\II{\sh{\L}{\si(r)}{}+1}{\sh{\L}{\si(r-1)}{}}
\bs\OT (\L),
}
(recall that $\OT ({\l^{(r-1,i)}})=\OT (\L)$),
such that
an integer in this set is skipped.
So
\equa{2-numbers}{
\#\bigl(\II{\sh{\L}{\si(r)}{}}{\sh{\L}{\si(r-1)}{}}
\bs\OT (\L)\bigr)\ge2.
}
In particular there is at least an element in the set (\ref{the-set}),
which means
\equa{equal-to-0}
{
\dd^\L_{\si(r),\si(r-1)}=0\mbox{ \ \ \ if \ }\si(r-1)-\si(r)=1
\mbox{ \ \ \ (cf.~(\ref{c'-st}))}.}
Note that $\sh{\L}{\si(r)}{}$ is in
$Q_{r-1}=
\II{\sh{\l}{r-1}{'}+1}{\sh{\L}{\si(r-1)}{}}
\bs \OT (\L)$
by (\ref{entry-r-1}), but not in the set (\ref{the-set}).
Thus we have in fact
skipped two integers in the set $Q_{r-1}$. Therefore
$$
\mbox{$\theta_{r-1}=\theta'_{r-1}-2$
\ \ \ \ if \ \ $\si(r-1)>\si(r).$}
$$
\vskip4pt
{\it Case 4: The general case.}

In general, when we apply $R_{\mm'_s,\cp\nn'_s}$ to $\l^{(s)}$ in order to obtain $\l^{(s-1)}$, for
each $t>s$ with $\si(t)<\si(s)$ the above arguments show that
we have to skip two integers in the set
$Q_s$. Therefore we have (\ref{theta-s}) and the lemma.
\end{proof}
\begin{lemma}\label{t-l-3}
The image of the map $(\ref{injection})$ is contained in $\fS^{\L,\l}$, i.e., if $\si=\si_\theta$
for some $\theta\in\Theta^\L_\l$, then
\equa{to-prove}
{
\dd^\L_{\si(t),\si(s)}=0\mbox{ \ \ \ if \ \ }\si(t)<\si(s)
\mbox{ \ \ for any \ \ }t>s.
}
\end{lemma}
\begin{proof}
Let $p$ be the number such that $\si(p)$ is minimal
among those $\si(u)$ with $u<t$ and $\si(u)>\si(t)$.
Then $p=s$ or $\si(p)<\si(s)$.
The definition (\ref{define-d-st}) means that the relation
$\dd^\L_{\si(t),\si(p)}=0$ implies the relation (\ref{to-prove}). Thus
it suffices to prove $\dd^\L_{\si(t),\si(p)}=0$.

For any $p'$ with $\si(t)\le \si(p')<\si(p)$, by the choice of
$p$, we have $p'\ge t$ and so $p'>p$.
Thus the arguments in the proof of Lemma \ref{t-l-2}
show that when we apply $R_{\mm'_p,\cp\nn'_p}$
to $\l^{(p)}$ in order to obtain $\l^{(p-1)}$,
we need to pass over the integer $\sh{\L}{\si(p')}{}$
and skip another integer for all such $p'$. Hence there are
at least $2(\si(p)-\si(t))$ integers in the
set
$\II{\sh{\L}{\si(t)}{}}{\sh{\L}{\si(p)}{}}
\bs\OT (\L)$ (cf.~(\ref{2-numbers})). This means that
there are at least
$\si(p)-\si(t)$ integers in the
set
$\II{\sh{\L}{\si(t)}{}}{\sh{\L}{\si(p)}{}}
\bs \Set(\L^\rho)$
(since there are exactly $\si(p)-\si(t)$ integers in
$\II{\sh{\L}{\si(t)}{}}{\sh{\L}{\si(p)}{}}$ which are
entries of atypical tuple $\tt_\L$), which implies that
$\cc^\L_{\si(t),\si(p)}=0$ by (\ref{c'-st}) and so
$\dd^\L_{\si(t),\si(p)}=0$
(cf.~(\ref{equal-to-0})).
\end{proof}
\begin{lemma}\label{t-l-4}
The map $(\ref{injection})$ is a bijection between $\Theta^\L_\l$ and $\fS^{\L,\l}$.
\end{lemma}
\begin{proof}
For any $\si\in \fS^{\L,\l}$, we define $\theta_s$ as in
(\ref{theta-s}).
We want to prove
\equa{to-p}
{\theta_s\ge0\mbox{ \ \ \ for \ \ }s=1,...,r.
}
Suppose $t>s$ such that $\si(t)<\si(s)$.
So
\equa{so1}
{
\dd^\L_{\si(t),\si(s)}=0.
}
Denote
\equa{use0}
{
X_{s,t}=\{p>s\,|\,\si(t)\le\si(p)<\si(s)\}\mbox{ \ \ and \ \ }
x_{s,t}=\#X_{s,t}.
}
First we prove by induction on $\si(s)-\si(t)$ that
\equa{use01}
{
\#\bigl(
\II{\sh{\L}{\si(t)}{}}{\sh{\L}{\si(s)}{}}
\bs \Set(\L^\rho)\bigr)\ge x_{s,t}.
}
If
$\si(s)-\si(t)=1$, then (\ref{use01}) follows from (\ref{so1})
and Definition \ref{c-related}
%,(\ref{c'-st}), (\ref{c'-st})$'$ and (\ref{define-d-st})
(note that obviously $\si(s)-\si(t)\ge x_{s,t}$).

In general set $p'=x_{s,t}$ and we write
the set $\{\si(p)\,|\,p\in X_{s,t}\}$ in ascending order:
\equa{si-p}
{
\{\si(p)\,|\,p\in X_{s,t}\}
=\{\si(t_1)<\si(t_2)<\cdots<\si(t_{p'})\},
}
where $t_1=t$ and $\si(t_{p'})<\si(s).$
Since $\dd^\L_{\si(t),\si(s)}=0$,
by Definition \ref{c-related}, %(\ref{c'-st}), (\ref{c'-st})$'$ and (\ref{define-d-st}),
there exists some number, denoted by $\si(i)$, lies
in between $\si(t)$ and $\si(s)$,
i.e.,
\equa{i.e.}{
\si(t)<\si(i)\le\si(s),
}
such that
$\cc^\L_{\si(t),\si(i)}=0$, that is,
\equa{use1}
{
\#\bigl(\II{\sh{\L}{\si(t)}{}}{\sh{\L}{\si(i)}{}}\bs \Set(\L^\rho)\bigr)\ge \si(i)-\si(t)
\mbox{ \ \ \ (cf.~(\ref{d(s,t)}))}.
}
Now we prove (\ref{use01}) by induction on $\si(s)-\si(t)$ in two cases.
\vskip6pt
{\it Case 1: Suppose $i\in X_{s,t}$, say $i=t_{p''}$ for some $1<p''\le p'$}.
Then
\equa{at-least-p'}
{\#\bigl(
\II{\sh{\L}{\si(t_{p''})}{}}{\sh{\L}{\si(s)}{}}
\bs \Set(\L^\rho)\bigr)\ge
x_{s,t_{p''}}
\ge
 p'-p''+1,
}
where the first inequality follows from
the inductive assumption that (\ref{use01}) holds for
$t_{p''}$
since $\si(s)-\si(t_{p''})
<\si(s)-\si(t)$, and the second %inequality follows
from
the fact that $t_{p''},t_{p''+1},...,t_{p'}\in$ $X_{s,t_{p''}}$.
Then (\ref{use01}) follows from (\ref{at-least-p'}) and
(\ref{use1}) (with $i$ replaced by $t_{p''}$) by
noting that $\si(t_{p''})-\si(t)\ge
p''-1$ (cf.~(\ref{si-p})) and that $p'=x_{s,t}$.

\vskip6pt
{\it Case 2: Suppose $i\notin X_{s,t}$.}
 This means that $i<s$.
Let $p''$ be the minimal integer with $1\le p''\le p'$
such that $\si(t_{p''})>\si(i)$.
Then $\si(s)-\si(t_{p''})<\si(s)-\si(t)$ and so
(\ref{at-least-p'}) holds again in this case by the inductive assumption. Furthermore
since $\si(i)-\si(t)<\si(s)-\si(t)$ (cf.~(\ref{i.e.})),
the inductive assumption also gives
\equa{t1-i}{
\#\bigl(
\II{\sh{\L}{\si(t)}{}}{\sh{\L}{\si(i)}{}}
\bs \Set(\L^\rho)
\bigr)\ge x_{i,t}\ge p''-1,
}
where the last inequality follows from the fact that
$$
t_1,...,t_{p''-1}\in \{p>i\,|\,\si(t)\le\si(p)\le\si(i)\}=X_{i,t},$$
(recall that $i<s$).
Now (\ref{use01}) follows from (\ref{t1-i}) and
(\ref{at-least-p'}). This completes the proof of (\ref{use01}).
\vskip6pt
Now set $p=\ell(\si,s)$ and write
\equa{w1}
{
\{t>s\,|\,\si(t)<\si(s)\}=
\{s_1,s_2,...,s_p\,|\,\si(s_1)<\si(s_2)<\cdots<\si(s_p)<\si(s)\}.
}
(Then the left-hand side of (\ref{w1}) is in fact the
set $X_{s,s_1}$, cf.~(\ref{si-p}).) Thus
(\ref{use01}) means that the set
\equa{w2}
{
\II{\sh{\L}{\si(s_1)}{}}{\sh{\L}{\si(s)}{}}
\bs \Set(\L^\rho)
}
has cardinality $\ge x_{s,s_1}= p$,
and
so the set
\equa{w3}
{
\II{\sh{\l}{s}{'}+1}{\sh{\L}{\si(s)}{}}
\bs \OT (\L)
}
has cardinality $\ge 2p,$
because $\sh{\l}{s}{'}
<\sh{\l}{s_i}{'}
\le \sh{\L}{\si(s_i)}{}
<\sh{\L}{\si(s)}{},$
i.e., we have $p$ more elements $\sh{\L}{\si(s_i)}{}$
with $i\in\II{1}{p}$ which are not in the first set (\ref{w2}) but in the second set
(\ref{w3}).

Then (\ref{Q-s}), (\ref{ell(si,s)}) and the fact that (\ref{w3}) has cardinality $\ge 2p$
show that
$\theta'_s\ge 2p$ and so $\theta_s\ge0$.
(In fact if $p>0$ then $\theta_s\ge1$ since in this case there exists at least one more
integer $\sh{\L}{\si(s)}{}$ which is in
$\II{\sh{\l}{s}{'}+1}{\sh{\L}{\si(s)}{}}
\bs \OT (\L).$) This proves (\ref{to-p}).

Now we define $\l^{(0)}$ as in (\ref{denote-gi}). Then the arguments
in the proof of Lemma \ref{t-l-2} show that
(\ref{g1-si-f}) holds, i.e., $(\l^{(0)})^+=\L$. Thus $\theta\in\Theta^\L_\l$ and $\si=\si_\theta$.
Therefore
$\theta\mapsto\si_\theta$ is a bijection
between $\Theta^\L_\l$ and $\fS^{\L,\l}$.
\end{proof}

\ni{\it Proof of Theorem \ref{theo7.1}.} Finally we return
to the proof of Theorem \ref{theo7.1}.
By (\ref{l-g-f}), (\ref{Fact2}), (\ref{ell-si}) and (\ref{Q-s}), we have
$$\mbox{$
|\theta|=\sum\limits_{s=1}^r\theta'_s
-\sum\limits_{s=1}^r\ell(\si,s)=\ell(\L,\l)-\ell(\si).
$}$$
Now (\ref{k-l-p-f}) follows from (\ref{ell-gf}).\hfill$\Box$
%\vskip8pt

\def\Ob{{\rm Ob\,}}\def\OL{\overline}
\def\E{{\mathcal{E}}}\def\B{{\mathcal{B}}}\def\deg#1{{\# #1}}
\subsection{A correspondence between $r$-fold atypical modules over $\gl_{m|n}$ and $\gl_{r|r}$ }
A subset $\B\subset \Z^{m|n}_+$ of dominant integral weights is called a {\it block
of $\Z^{m|n}_+$ for $\gl_{m|n}$}
if it is a maximal subset such that for any two weights $\L,\l\in \B$, there exist weights $\l^1,\l^2,...,
\l^k\in\B$ with $\L=\l^1,\,\l^k=\l$ such that the extension group
${\rm Ext}^1(V(\l^i),V(\l^{i+1}))\ne0$ for $i=1,...,k-1$.
Then $\Z^{m|n}_+$ is divided into a disjoint union of blocks.
Lemma 1.12 in \cite{Se96} says that for any $\L,\l\in\B$, one has
$\#\L=\#\l$,
which is called the {\it degree of atypical type of $\B$}, and denoted by
$\deg\B$.
Let $\E^{m|n}$
be the category of finite dimensional $\gl_{m|n}$-modules.
A dominant weight $\l$ is called a {\it primitive weight of
a module $V$} if it is the highest weight of a composition factor
of $V$.\def\PW{{\rm Prim}}
Denote by $\PW(V)$ the set of primitive weights of $V$.
For $\L,\l\in\Z_+^{m|n}$, we denote by
$a_{\L,\l}=[\OL V(\L):V(\l)]$ the multiplicity of the composition
factor $V(\l)$ in the Kac-module $\OL V(\L)$. It was proved in \cite{B} that
\equa{a-l-mu}
{
a_{\L,\l}\le1,\mbox{ \ and the matrix $(a_{\L,\l})$ is the inverse of
the matrix $(K_{\L,\l}(-1))$},
}
where the matrices were defined with respect to some total order of
weights compatible with the partial order ``$\soe$''.

%For any block $\B$, we denote by $\E^{m|n}_\B$ the full subcategory of
%$\E^{m|n}$ such that every primitive weight of any module in $\E^{m|n}_\B$ belongs to $\B$.
%%
%%Two subcategories $\E_1,\E_2$ of $\E^{m|n},\E^{m'|n'}$ are {\it equivalent} if there
%%is a one to one correspondence between $\E_1$ and $\E_2$ such that ...

An application of Theorem \ref{theo7.1} is the following.
\begin{theorem}\label{theo7.1'}
Let $\B$ be a block of $\Z_+^{m|n}$ for $\gl_{m|n}$ with
$\deg\B=r$. Let $\B'$ be the unique block
of $\Z_+^{r|r}$ for $\gl_{r|r}$ with $\deg{\B'}=r$.
\begin{enumerate}
\item\label{theo7.1'-1}
There exists a bijection
\equa{sends}
{
\phi:\B\to\B',\ \ \ \
\phi(\L)=(h'(\L)\,|\,h(\L)),
}
where $h(\L)$ is the height vector of $\L$ defined by $(\ref{height})$ and
$h'(\L)\!=\!(h_r(\L),...,h_1(\L))$. % is the reverse vector of $h(\L)$,
\item \label{theo7.1'-2}
The Kazhdan-Lusztig polynomials $K_{\L,\l}(q)$ of $\,\gl_{m|n}$ and
$K_{\phi(\L),\phi(\l)}(q)$ of $\,\gl_{r|r}$ coincide, that is,
\equa{K-L-equal}
{
K_{\L,\l}(q)=K_{\phi(\L),\phi(\l)}(q)
}
for $\L,\l\in\B.$
\item \label{theo7.1'-3}
%%Thus by \cite [{\rm Theorem 4.51}]{B}, the
%The category $\E^{m|n}_\B$ is equivalent to
%the category $\E^{r|r}_{\B'}$; in particular all subcategories $\E^{m|n}_\B$ of
%$\E^{m|n}$ with $\deg\B=r$ are equivalent.
Under the mapping $(\ref{sends})$,
the set $\PW(\OL V(\L))$ of primitive weights of the $\gl_{m|n}$-Kac-module $\OL V(\L)$ for $\L\in\B$ is
in one to one correspondence to the set  $\PW(\OL V(\phi(\L)))$
of primitive weights of the $\gl_{r|r}$-Kac-module $\OL V(\phi(\L))$, namely
\equa{kac-pri}
{
\PW(\OL V(\phi(\L)))=\{\phi(\l)\,|\,\l\in\PW(\OL V(\L))\}.
}
\end{enumerate}
\end{theorem}
\ni{\it Proof.~}~(\ref{theo7.1'-1})
From Theorem \ref{theo7.1} (cf.~(\ref{prec})) we can deduce
that any two dominant weights are in the same block if and only if their typical tuples
are equal (thus blocks of $\Z_+^{m|n}$ are in one to one correspondence with the
typical tuples). Therefore for any block $\B$ of $\Z_+^{m|n}$ with $\deg\B=r$,
we can denote the typical tuple of any weight in $\B$ by
\def\tot{{\rm tt}}\equa{denote-typical}
{
\tot_\B=(\tot_1,...,\tot_{m-r}\,|\,\tot_{\cp r+1},...,\tot_{\cp n}),
}
where $\tot_1>...>\tot_{m-r},\,\tot_{\cp r+1}<...<\tot_{\cp n}$ are all distinct.

First we prove that the map $\phi$ is injective, which is equivalent to proving that
for any two dominant weights $\L,\l\in\B$ with $h(\L)=h(\l)$, one has the equality of
the atypical tuples $\tt_\L=\tt_\l$.
Thus assume that
$\L^\rho_{\cp\nn_s}\ne\l^\rho_{\cp\nn'_s}$ for some $s$. Say,
$\L^\rho_{\cp\nn_s}>\l^\rho_{\cp\nn'_s}$. Then the right-hand side of (\ref{prop-1-1})
with $t=s$ is zero, but the left-hand side is not since
$\L^\rho_{\cp\nn_s}\in\II{\l^\rho_{\cp\nn'_s}+1}{\L^\rho_{\cp\nn_s}}\bs\Set(\tot_\B)$.

Next we prove that the map $\phi$ is surjective, which is equivalent to proving that
for any lexical $r$-tuple $a=(a_1,...,a_r)\in\Z^r$, there exists $\L\in\B$ such that
$h(\L)=a$. Fix a lexical $r$-tuple $b=(b_1,...,b_r)\in\Z^r$ such that
$b_r<\min\{\tot_{m-r},\tot_{\cp r+1},a_1\}-r$. Let
\equa{let-mu}
{
\l=-\rho+(\tot_1,...,\tot_{m-r},b_r+r,...,b_1+1\,|\,
b_1+1,...,b_r+r,\tot_{\cp r+1},...,\tot_{\cp n})\in\Z^{m|n}.
}
Then $\l$ is dominant such that $\ot_\l=\tot_\B$ (thus $\l\in\B$)
and $h(\l)=b$. Take $t=s$ and
$h_s(\L)=a_s$ in (\ref{prop-1-1}), then (\ref{prop-1-1}) uniquely determines
a number $\L^\rho_{\cp\nn_s}\notin\Set(\tot_\B)$.  This uniquely determines
a weight $\L\in\B$. To see $h(\L)=a$, use (\ref{prop-1-1}) again.
Thus $\phi$ is a bijection.

In fact the atypical tuple $\tt_\L$  determines the height vector $h(\L)$
through the following formula:
\equa{a-h(L)}
{
h_s(\L)+s=\L^\rho_{\cp\nn_s}-\#\{\mbox{typical entries which are smaller than $\L^\rho_{\cp\nn_s}$}\}.
}
Conversely  the height vector $h(\L)$  determines  the atypical tuple $\tt_\L$ as follows:
First set $\tt_\L=h(\L)$, and denote $\tt_\L$ as $\tt_\L=(a_1,...,a_r)$.
Label the typical entries of $\ot_\L$ in ascending order:
$x_1<x_2<...<$ $x_{m+n-r}$. For each $i=1,2,...,m+n-r$, if $a_t\ge x_i$ then replace $a_t$ by $a_t+1$
for all $t=1,2,...,r$.

(\ref{theo7.1'-2})
Note that in $\gl_{r|r}$, we have $\tt_{\phi(\L)}=(\phi(\L))^{\rho_{r|r}}$
(cf.~notation (\ref{l-rho})),
where
$\rho_{r|r}=$ $(r,...,1\,|\,1,...,r)$,
and the typical tuple $\ot_{\phi(\L)}$ is empty. Also we have the equality of the height vectors:
$h(\L)=h(\phi(\L))$. Thus
by Definition \ref{c-related}, we have
$\cc^{\L}_{s,t}=\cc^{\phi(\L)}_{s,t}$ for all $s,t\in\II{1}{r}$.
Let $\l$ be another dominant weight with $\l\soe\L$.
By (\ref{prop-1-1}) and (\ref{where}), we have
$b^{\L,\l}=b^{\phi(\L),\phi(\l)}$.  This gives
(\ref{K-L-equal}) by
(\ref{F-x'}) and (\ref{k-l-p-f}).

(\ref{theo7.1'-3})
Using (\ref{a-l-mu}) and (\ref{K-L-equal}), we have
\\[6pt]\hs{3ex}
$
\PW(\OL V(\phi(\L)))=\{\l'\,|\,a_{\phi(\L),\l'}=1\}
$\\[6pt]\hs{3ex}$\phantom{\PW(\OL V(\phi(\L)))}
=\{\phi(\l)\,|\,a_{\phi(\L),\phi(\l)}=1\}
=\{\phi(\l)\,|\,a_{\L,\l}=1\}
=\phi(\PW(\OL V(\L)))
.
$\hfill$\Box$

%%%%%%%%%%%%%%%%%%%%%%%%%%%%%%%%%%%%%%%%%%%%%%%%%%%%%%%%%%%%%%%%
%%%%%%%%%%%%%%%%%%%%%%%%%%%%%%%%%%%%%%%%%%%%%%%%%%%%%%%%%%%%%%%%
%
\section{Character formulae}
\label{char-for}
%%%%%%%%%%%%%%%%%%%%%%%%%%%%%%%%%%%%%%%%%%%%%%%%%%%%%%%%%%%%%%%%
\def\L{\Lambda}\def\l{\lambda}
As mentioned in the introduction, this section contains three main results: the proof of the conjecture due
to van der Jeugt et al, the construction of a Kac-Weyl type character formula, and
the derivation of a dimension formula.

We shall continue to use notations in the previous sections. Moreover, we
define
\equa{m(f)-g}
{
m(\L)_\l=\#\fS^{\L,\l}}
to be the cardinality of the set $\fS^{\L,\l}$ (cf.~(\ref{s-f-r})).
Then $m(\L)_\l=\fS^{\L,\l}(1)$ (cf.~(\ref{l-q-function})).

For any weight $\l\in\PP$, we define
what is called the {\it Kac-character of $\l$}:
\equa{kac-ch}
{
\CHIK(\l)=\frac{L_1}{L_0}\mbox{$\sum\limits_{w\in W}$}\es(w)e^{w(\l+\rho)}.
}
Namely, $\CHIK(\l)$ is defined by the right-hand side of
(\ref{typical-char}). % with $\L$ replaced by $\l$.
Thus it is the character of the Kac-module $\overline V(\l)$
when $\l$ is dominant.

\subsection{Proof of the conjecture of van der Jeugt et al}
As an immediate consequence of Theorem \ref{theo7.1}, we have
\begin{theorem}\label{ours}
The formal character $\ch  V(\L)$ of the finite dimensional
irreducible $\fg$-module $V(\L)$
is given by
\equa{ch-v}{
\ch V(\L)
%=\mbox{$\sum\limits_{\l\in \PP _+}$}K_{\L,\l}(-1)\CHIK(\l)
=
\mbox{$\sum\limits_{\l\in \PP _+:\:\l\soe\L}$}(-1)^{\ell(\L,\l)}
m(\L)_\l%\mbox{$\sum\limits_{\si\in \fS^{\L,\l}}$}
\CHIK(\l),
}
where $\PP_+$ is the set of dominant integral weights, the length $\ell(\L,\l)$
is defined in $(\ref{l-g-f})$, and the partial order ``$\soe$''
is defined in $(\ref{g-prec-f})$.
\end{theorem}
\noindent{\it Proof.} %\begin{proof}
This follows from  (\ref{k-l-p-f}), (\ref{m(f)-g}) and
\cite[Lemma 3.4]{Se98}, which states
$$\mbox{$
\ch V(\L)=\sum\limits_{\l\in \PP _+}K_{\L,\l}(-1)\CHIK(\l).
$}
\eqno\Box
$$
%\end{proof}
\vskip4pt
One can re-write (\ref{ch-v}) to obtain the conjecture of van der Juegt et al. To state it,
we need to introduce the following notations.

Let $\l$  as in (\ref{(1)}) be an $r$-fold atypical weight (not necessarily dominant)
with atypical roots ordered as in (\ref{aty-roots}): $\g_1<...<\g_r$.
We define the {\it normal cone with vertex $\l$}:
\equa{(13)}{
\begin{array}{c}
\mbox{$
\NC\l=\{\l-\sum\limits_{s=1}^r i_s\g_s\,|\,i_s\ge0\}$}.
\end{array}
}
We also define $\TC\l$,
which was referred to as the {\it truncated cone with vertex $\l$} in \cite{VHKT},
to be the subset of $\NC\l$
consisting of weights $\mu$ such that
the $s$-th entry of the atypical tuple $\tt_\mu$ (cf.~(\ref{aty-root-L}))
is smaller than or equal to
the $t$-th entry of $\tt_\mu$
when the atypical roots $\g_s,\g_t$ of $\l$ (not $\mu$)
are strongly $c$-related
for $s<t$. Namely (recall the $\rho$-translated notation in (\ref{l-rho}))
\equa{cf1}{
\begin{array}{c}
\mbox{$
\TC\l=\{\mu\in\NC\l\,|\,
\mu^\rho_{\cp \nn_s}\le
\mu^\rho_{\cp \nn_t}\mbox{ \ if \ }\dd^\l_{s,t}=1\mbox{ \ for \ }s<t\}.
$}
%\\[4pt]
\end{array}
}
For $\l=\L-\sum_{s=1}^r i_s\g_s
\in\NC\L$, we denote
\begin{eqnarray}
\label{relative}
|\L-\l|=\mbox{$\sum\limits_{s=1}^r$} i_s
\mbox{ \ \ (called the {\it relative level of $\l$})}.
\end{eqnarray}

As an application of Theorem \ref{ours} we prove the
following character formula which was a conjecture
put forward by van der Jeugt, Hughes, King and
Thierry-Mieg in \cite{VHKT} as the result of in depth research carried out by the authors over several years time.
%%%%%%%%%%%%%%%%%%%%%%%%%%%%%%%%%%%%%%%%%%%%%%%%%%%%%%%%
\begin{theorem}\label{Conj-for}
%%%%%%%%%%%%%%%%%%%%%%%%%%%%%%%%%%%%%%%%%%%%%%%%%%%%%%%%
\begin{eqnarray}
\label{formula}&\!\!\!\!\!&\!\!\!\!\!\!\!\!
\mbox{$
\ch V(\L)
=\sum\limits_{\l\in \TC\L}(-1)^{|\L-\l|}\CHIK(\l).
$}
\end{eqnarray}
\end{theorem}
\begin{proof}
First we remark that for %any $w\in W$ and
any weight $\l\in \PP $, we have
(cf.~(\ref{kac-ch}))
\equa{vanishing}
{
\CHIK(\l)=0\mbox{ \ \ \ if $\l$ is vanishing (i.e., not regular)}.
}
Note that for any weight $\l$ in the truncated cone $\TC\L$ or in
the set
\equa{R-Lambda}
{R_\L:=\{\mu\in \PP _+\,|\,\mu\soe\L\}
\mbox{ \ \ (cf.~the right-hand side of (\ref{ch-v}))},
}
we have the equality of the typical tuples:
\equa{t-equal}{
\ot_{\l}=\ot_{\L}
\mbox{ \ \ \ (cf.~(\ref{ot-f}))}.
}
Thus $\l$ is uniquely
determined by the atypical tuple $\tt_{\l}$ (cf.~(\ref{aty-root-L})).
Clearly for any $\mu\in \TC\L$, it is either vanishing
(so $\CHIK(\mu)=0)$
or is $W$-conjugate under the dot action
to some unique $\l\in R_\L$, where the later fact is equivalent to the existence of
a unique $\si\in \Sr_r$ satisfying
$\tt_{\si\cdot\mu}=\tt_{\l}$
(cf.~(\ref{(4)}) and (\ref{Action-of-s-r})).
If $\si(t)<\si(s)$ for some $1\le s<t\le r$. Then
$$
\mu^\rho_{\cp\nn_{\si(t)}}
=\l^\rho_{\cp\nn_t}
>\l^\rho_{\cp\nn_s}
=
\mu^\rho_{\cp\nn_{\si(s)}}.
$$
Thus $\dd^\L_{\si(t),\si(s)}=0$ by (\ref{cf1}).
Therefore $\si\in S^\L$ by (\ref{s-f-r}).
Obviously $\tt_{\l}=
\tt_{\si\cdot\mu}\le$
$\tt_{\si\cdot\L}$ (cf.~(\ref{partial-order-on-Z-r})),
i.e., $\si\in \fS^{\L,\l}$ by definitions (\ref{s-f-g}) and (\ref{g-prec-f})
and the fact that $\tt_{\si\cdot\L}=\tt_\L=$ $\tt_\l$ (cf.~(\ref{t-equal})).
Conversely, for any $\l\in R_\L$ and $\si\in \fS^{\L,\l}$,
there corresponds to a unique
$\mu\in \TC\L$ such that $\tt_{\si\cdot\mu}=\tt_{\l}$.

For any $\l\in R_\L$,
%let $M(\l)
%=\{\mu\in \TC\L\,|\,\mu\mbox{ is $W$-conjugate under the dot
%action to }\l\},
%=\TC\L\cap W\cdot\l,$ and
let $\tilde\l\in \TC\L\cap W\cdot\l$ be the (unique)
lexical weight in the sense of
Definition \ref{defi-lexical}. Then
%each $\mu\in M(\l)$ is conjugate under the dot action to $\tilde\l$
%by an element of $\Sr_r$, thus by an element of $W$ with even parity
%(cf.~(\ref{(4)}) and (\ref{Action-of-s-r})). Therefore
\equa{chi-equal}
{\CHIK(\mu)=\CHIK(\tilde\l)
\mbox{ \ \ \ for \ }\mu\in  \TC\L\cap W\cdot\l,
}
since elements of $\Sr_r$ correspond to elements of $W$ with even parity
(cf.~(\ref{(4)}) and (\ref{Action-of-s-r})).

The above arguments have in fact shown that
the right-hand side of (\ref{formula}) is equal to
\equa{fur1}
{
\mbox{$\sum\limits_{\l\in \PP _+:\:\l\soe \L}$}(-1)^{|\L-\l|}
\mbox{$\sum\limits_{\si\in \fS^{\L,\l}}$}
\CHIK(\tilde\l)=
\mbox{$\sum\limits_{\l\in \PP _+:\:\l\soe \L}$}(-1)^{|\L-\l|}
m(\L)_\l
\CHIK(\tilde\l).
}
By (\ref{ch-v}), what remains to prove is the following:
if $w(\tilde\l^\rho)=\l^\rho $ for $w\in W,\,\l\in R_\L$, then
\equa{to-proof}
{
(-1)^{|\L-\l|}=(-1)^{\ell(\L,\l)}\es(w).
}
Note that in order to obtain $\l^\rho $ from
$\tilde\l^\rho$ by moving
all entries
$\tilde\l^\rho_{\mm_s}$, $\tilde\l^\rho_{\cp\nn_s}$ for $s=1,...,r$
to suitable positions
step by step, each time exchanging nearest neighbor entries only,
the total number of movements is $N=\sum_{s=1}^r N_s$, where
$$
N_s=\#(\II{\tilde\l^\rho_{\cp\nn_s}+1}{\L^\rho_{\cp\nn_s}}
\cap\OT(\L))
\mbox{ \ \ (cf.~notation (\ref{set-R(f)}))}.
$$
To see this, note that $\tilde\l$ is obtained from $\L$ by subtracting
some atypical roots (cf.~(\ref{(13)})).
Thus the entries $\tilde\l^\rho_{\mm_s},\,
\tilde\l^\rho_{\cp\nn_s}$ of $\tilde\l^\rho$, which seat at the
positions
$\mm_s,\cp\nn_s$ of $\L$ should be moved to appropriate positions
in order to make the resultant weight $\l^\rho$
dominant. The
number of steps needed for these two entries is obviously $N_s$.

From (\ref{l-g-f}) and (\ref{l-g-f'}),
we have $|\L-\l|=\sum_{s=1}^r(\L_{\cp\nn_s}-\l_{\cp\nn_s})=
\ell(\L,\l)+N$. But $\es(w)=(-1)^N$, we obtain (\ref{to-proof}).
\end{proof}

\subsection{Definitions of $\l\lm$ and $C_r$}
Our purpose is to re-write (\ref{formula}) into a finite sum.
Since the sum over the truncated cone $\TC\L$ in (\ref{formula})
is difficult to compute, we want to change
this sum into several sums over some normal cones $\NC\mu$
by making use of the fact that
the Kac-character $\CHIK(\l)$ is
$\Sr_r$-invariant under the dot action (cf.~(\ref{(4)})
and (\ref{Action-of-s-r})), as a
%Then the
sum over a cone $\NC\mu$ is easy to compute.

This will be done in two steps.

First we need to introduce more notations.
Define another partial order ``$\le$'' on $\NC\L$ such that for  $\l,\mu\in\NC\L$,
\equa{part}
{\mu\le\l  \ \ \ \Lra \ \ \
\mbox{every entry of $\mu$
$\ \ \le\ \ $
the corresponding entry of $\l$}.
}
%In the rest of the paper we shall use this partial order ``$\le$'' on $\NC\L$ instead
%of the previous one ``$\soe$'' defined in (\ref{g-prec-f}).
\begin{notation}\label{Notation1}
\rm
For $\l\in\NC\L$, denote
by $\l\lm\in\NC\L$
the maximal lexical weight (cf.~(\ref{lexical})) which is $\le\l$, namely,
\begin{eqnarray}
\label{l-m}
\l\lm&\!\!\!\!=\!\!\!\!&
\max\{\mu\in\NC\L\,|\,\,\mu\le\l,\,
\mbox{ \ and \ }\mu\mbox{ is lexical}\}.
\end{eqnarray}
Thus we have the equality of
typical tuples: $\ot_{\l\lm}=\ot_\l$ (cf.~(\ref{ot-f})) and the
entries of atypical tuple $\tt_{\l\lm}$ (cf.~(\ref{denote-t-f})) are defined by
\equan{l-m'}
{
\shb{\l\lm}{s}{}=\min\{\sh{\l}{t}{}\,|\,s\le t\le r\}
\mbox{ \ \ \ for \ \ }s=1,...,r.
%\eqno \PED
}
%where $\tt_{\l\lm}=(\sh{\l\lm}{1}{},...,\sh{\l\lm}{r}{})$ and
%$\tt_\l=(\sh{\l}{1}{},...,\sh{\l}{r}{})$ (cf.~notation (\ref{s-h})).
%\PED
\end{notation}

Denote by $\HC\l$,
called the {\it lexical cone with vertex $\l$},
%to be
the subset of the truncated cone $\TC\l$ (cf.~(\ref{cf1})) consisting of
lexical weights (cf.~(\ref{lexical})),
namely
\equa{(14)}{
\begin{array}{l}
\mbox{$
\HC\l=\{\mu\in\TC\l\,|\,\,
%\mu^\rho_{\cp \nn_s}\le \mu^\rho_{\cp \nn_t}\mbox{ for }s<t
\mu\mbox{ is lexical}
\}.
$}
%\\[4pt]
\end{array}
}

%\vskip4pt
%{\bf Step 1.~}~
Our
first step is to
change the sum over $\TC\L$ in (\ref{formula}) to several
\mbox{sums} over some lexical cones $\HC{(\si\cdot\L)\lm}$
(see Proposition \ref{prop1}).

%\vskip4pt
The proof of Theorem \ref{theo7.1} (cf.~the arguments of proving
(\ref{fur1})) and (\ref{vanishing}) show that
(\ref{formula}) can be re-written as
\equa{(15)}
{\mbox{$
\ch V(\L)=
\sum\limits_{\si\in S^\L}\ \
\sum\limits_{\l\in \HC\L:\:\l{\ssc\,}\soe{\ssc\,}
\si\cdot\L,{\sc\,}\l\mbox{ regular}}\
(-1)^{|\L-\l|}\CHIK(\l).
$}}
The definition of $\l\lm$ in (\ref{l-m}) shows that we have the equality of
the following two sets of regular weights:
\equa{(16)}{
\{\l\in\HC\L\,|\,\l\soe\si\cdot\L,\,\l\mbox{ is regular}\}
=\{\l\in\HC{(\si\cdot\L)\lm}\,|\,\l\mbox{ is regular}\}.
}
Thus (\ref{(15)}) leads to the following.
\begin{proposition}\label{prop1}
\equa{(17)}
{\mbox{$
\ch V(\L)=
\sum\limits_{\si\in S^\L}\ \
\sum\limits_{\l\in \HC{(\si\cdot\L)\lm}}\
(-1)^{|\L-\l|}\CHIK(\l).
$}}
\end{proposition}
%which is what we want in the first step.

The second step is to change
the sum over the lexical cone
$\HC{(\si\cdot\L)\lm}$ in (\ref{(17)}) to
several sums over the normal cones
$\NC{(\tau\cdot(\si\cdot\L)\lm)\lm}$.
This will be achieved by Lemma \ref{lemm7.2}
below. To state lemma, we need
some further notations.

\begin{notation}\label{Notation3}
\rm
Denote by $C_r$ the subset of $\Sr_r$ consisting of
permutations $\tau$ which can be written as a product of
{\it cyclic permutations}
of the form
\equa{C-r}
{{\sc}
\tau=(1,2,...,i_1)(i_1+1,i_1+2,...,i_1+i_2)\cdots
(i_1+...+i_{t-1}+1,i_1+...+i_{t-1}+2,...,r),
%i_1\!+\!...\!+\!i_{t}),
}
where $i_1,...,i_t$ are positive integers such that
$\sum_{s=1}^ti_s=r$
(namely, $(i_1,...,i_t)$ is a {\it composition of $r$}).
Associated to $\tau$, there is the multi-nomial coefficient
\equa{length}
{
\left({\,}^{^{\dis r}}_{_{\dis \tau}}{\ssc\,}\right)
=\frac{r!}{i_1!\cdots i_t!}.
%\vs{-8pt}
}
%\vs{-6pt}\PED
\end{notation}
\def\hM{{\mathcal C}^{\rm Half}}
\def\hA{{\mathcal N}}
\def\hH{{\mathcal L}}
\def\hHH{{\mathcal H}}
\subsection{A technical lemma}
The following technical lemma is crucial in obtaining our character
formula in Theorem \ref{theo7.2}.
\begin{lemma}\label{lemm7.2}
Let $\l\in\HC\L$. We have
\equa{lemm7-2-1}
{\mbox{$
\sum\limits_{\mu\in\HC\l}(-1)^{|\L-\mu|}\CHIK(\mu)
={\dis\frac{1}{r!}}\sum\limits_{\tau \in C_r}
\left({\,}^{^{\dis r}}_{_{\dis \tau}}{\ssc\,}\right)
(-1)^{\ell(\tau )}
\sum\limits_{\mu\in\NC{(\tau\cdot\l)\lm}}(-1)^{|\L-\mu|}\CHIK(\mu),
$}}
where $\ell(\tau)$ is the length of $\tau$ $(cf.~(\ref{ell-si}))$, namely
$\ell(\tau)=\sum_{s=1}^t (i_s-1)$ for $\tau$ in $(\ref{C-r})$.
%Thus $(\ref{lemm7-2-1})$ is the inverse formula
%of $(\ref{chi-l})$.
\end{lemma}
\begin{proof}
For convenience, we denote by $C'_r$ the set of compositions of $r$, and denote
\equa{denote-tau-i}
{\tau_\bi= \
\tau \mbox{ which is defined by (\ref{C-r}) for \ }\bi=(i_1,...,i_t)\in C'_r.
}
We also denote (\ref{length}) by
$\left({\,}^{\sc r}_{\sc \bi}{\ssc\,}\right)$.
For any subset $S$ of $\NC\L$,
denote by $S\reg$ the regular weights of $S$.
We define
\equa{define-cC}
{
\chiP{S}:=\mbox{$\sum\limits_{\mu\in S}$}(-1)^{|\L-\mu|}\CHIK(\mu).
}
Thus
\equa{define-cC'}
{\chiP{S}=\chiP{S\reg}.
}
Note that any element $\mu\in\NC\L$
is uniquely determined by
the atypical tuple $\tt_{\mu}$ (cf.~(\ref{denote-t-f})), thus
there is a bijection
\equa{bijection-of}
{\mbox{$
\varphi:\L+\sum\limits_{s=1}^r\Z\g_s\to\Z^r,\ \ \
\varphi(\mu)=\tt_\mu=(\mu^\rho_{\cp\nn_1},...,\mu^\rho_{\cp\nn_r}).
$}}
For simplicity, we
use $g=\varphi(\l)=(g_1,...,g_r)$ to represent $\l$ (and transfer all
terminologies
to $g$).
Then the normal cone $\NC\l$ defined in (\ref{(13)}) and
the lexical cone $\HC\l$ defined in (\ref{(14)})
correspond to the
following sets respectively:
\begin{eqnarray}
\label{sets}
\NC g&\!\!\!\!=\!\!\!\!&
\{x=(x_1,...,x_r)\in\Z^r\,|\,x_s\le g_s,\,s=1,...,r\},
\nonumber\\
\HC g&\!\!\!\!=\!\!\!\!&
\{x\in\Z^r\,|\,x_1\le x_2\le\cdots\le x_s\mbox{ and }x_s\le g_s,\,s=1,...,r\}.
\end{eqnarray}
We define %what we call
the {\it half-lexical cone with vertex $\l$}:
\equan{denote-M'}
{
\hM_g=\{x\in\Z^r\,|\,x_2\le\cdots\le x_s\mbox{ and }x_s\le g_s,\,s=1,...,r\},
}
(i.e., we relax the condition of $x_1\le x_2$, cf.~(\ref{(14)}) and (\ref{sets})).
Define
\begin{eqnarray}
\label{denote-M}
&&\!\!\!\!\!\!\!\!\!\!\!\!\!\!\!\!\!\!
\hA(g_1,...,g_r)=\chiP{\NC g},\ \ \ \hH(g_1,...,g_r)=\chiP{\HC g},
\\
\nonumber
&&\!\!\!\!\!\!\!\!\!\!\!\!\!\!\!\!\!\!
\hHH_{g_1}(g_2,...,g_r)=\chiP{\hM_g},
\end{eqnarray}
which are the sign sums of Kac-characters over the normal cone $\NC g$,
the lexical cone $\HC g$ and
the half-lexical cone $\hM_g$ respectively (cf.~(\ref{define-cC})).
Then (\ref{lemm7-2-1}) is equivalent to
\equa{p2.1}
{\mbox{$
\hH(g_1,...,g_r)={\dis\frac{1}{r!}}\sum\limits_{\bi \in C'_r}
\left({\,}^{^{\dis r}}_{_{\dis \bi}}{\ssc\,}\right)
(-1)^{\ell(\bi )}\hA(g_\bi ),
$}}
where $g_\bi $ is defined by
\equa{l-composition}
{
g_\bi =\bigl({\sc\,}
\stackrel{i_1}{\overbrace{g_1,...,g_1}},{\sc\,}
\stackrel{i_2}{\overbrace{g_{i_1+1},...,g_{i_1+1}}},{\sc\,}...,{\sc\,}
\stackrel{i_t}{\overbrace{g_{\sum_{s=1}^{t-1}i_s+1},
{\sc\,}...,{\sc\,}g_{\sum_{s=1}^{t-1}i_s+1}}}
{\sc\,}\bigr)
}
for $\bi=(i_1,...,i_t)$.
This is because from the definition (\ref{l-m}), one can easily check
\equa{g-l-bi}
{g_\bi=\varphi((\tau_\bi\cdot\l)\lm) \mbox{ \ \ \
(cf.~(\ref{denote-tau-i}) and (\ref{bijection-of})).}
}
Denote
$$
\begin{array}{l}
g^{(s)}=(g_1,...,g_1,g_{s+1},...,g_r)=g|_{g_2=...=g_s=g_1}\in\Z^r,\\[4pt]
g'^{(s)}=(g_1,...,g_1,g_{s+1},...,g_r)\in\Z^{r-1},
\end{array}
$$
($g'^{(s)}$ is obtained from $g^{(s)}$ by deleting the first entry)
for $s=1,...,r$.
Note that
\equa{p2-2}{\mbox{$
\hHH_{g_1}(g'^{(1)})=\sum\limits_{s=1}^r \hH(g^{(s)}).
$}}
To see this, observe that
for any regular weight $x=(x_1,...,x_r)\in \hM_g$
(cf.~(\ref{define-cC'})),
there is a unique $s\in\II{1}{r}$ such that
\equa{con-x}
{x_2<...<x_s<x_1<x_{s+1}<...<x_r.
}
So $x$ is $\Sr_r$-conjugate
(thus under the inverse mapping $\varphi^{-1}$, it is $\Sr_r$-conjugate under the dot
action, cf.~(\ref{(4)}), (\ref{Action-of-s-r})) to
$$(x_2,...,x_s,x_1,x_{s+1},...,x_r)\in\HC{g^{(s)}}
\mbox{ \ \ since $x_1<g_1$.}
$$
Conversely every regular weight of $\HC{g^{(s)}}$ is $\Sr_r$-conjugate
under the dot action to
a unique weight $x$ of $\hM_g$ satisfying (\ref{con-x}). This proves (\ref{p2-2}).

By (\ref{p2-2}), we obtain
\equa{p2-3}{\mbox{$
\hH(g)=\hH(g^{(1)})=
\hHH_{g_1}(g'^{(1)})-\sum\limits_{s=2}^r (-1)^{s-1}{\dis\frac{1}{s(s-1)}}\hHH_{g_1}(g'^{(s)}).
$}}
Since the first variable $x_1$ in $\hM_g$ does not relate to any other variable,
when $g_1$ is fixed, $\hHH_{g_1}(g'^{(s)})$ is in fact $\hH(g'^{(s)})$
with respect to the $r-1$ variables $g_2,...,g_r$.
By the inductive assumption on $r$ that (\ref{p2.1}) holds for $r-1$, we have
\equa{p2.r-1}
{\mbox{$
\hHH_{g_1}(g'^{(s)})={\dis\frac{1}{(r-1)!}}\sum\limits_{{\bj }\in C'_{r-1}}
\left({\,}^{^{\dis r-1}}_{_{\dis \ \ \bj}}{\ssc\,}\right)
(-1)^{\ell({\bj })}\hA(g^{(s)}_{(1,{\bj })}),
\ s=1,...,r,
$}}
where
\equan{1-j}{
(1,\bj)=
\biggl\{\begin{array}{ll}
(1,j_1,j_2,...)&\mbox{if \ }s=1,
\\[4pt]
(j_1+1,j_2,...)&\mbox{otherwise},
\end{array}
\ \ \ \ \ \mbox{for \ \ }\bj=(j_1,j_2,...)\in C'_{r-1}.
}
Thus $(1,{\bj })$ is a composition of $r$.
For  $\bi \in C'_r$,
$g^{(s)}_\bi$
is defined by (\ref{l-composition}) with $g_2,...,g_s$ being set to $g_1$.
Thus each $g^{(s)}_{(1,{\bj })}$ has the form $g_\bi$
for some $\bi \in C'_r$, and
\equan{equal}
{
g^{(s)}_{(1,{\bj })}\!=\!g_\bi  \Lra
\biggl\{
\begin{array}{l}
i_1{\sc\!}={\sc\!}1,\,{\bj }{\sc\!}={\sc\!}(i_2,i_3,...)
\mbox{ if }s{\sc\!}={\sc\!}1,\mbox{ \ or}\\[4pt]
i_1{\sc\!}\ge{\sc\!}s,\,
1{\sc\!}+{\sc\!}j_1{\sc\!}+...{\sc\!}+{\sc\!}j_s
{\sc\!}={\sc\!}i_1\mbox{ and }(j_{s+1},j_{s+2},...){\sc\!}={\sc\!}(i_2,i_3,...)
\mbox{ otherwise}.
\end{array}
}
From this one can prove by induction on $r$ and $s$ that
the coefficient of $\hA(g_\bi )$ in (\ref{p2.r-1}) is
$\frac{1}{i_2!\cdots i_t!}b_{s,i_1}$, where $b_{0,i_1}=0$ and
$$
b_{s,i_1}=\frac{1}{i_1!}-\frac{1}{1!(i_1-1)!}+\frac{1}{2!(i_1-2)!}-...+
(-1)^{i_1-s}\frac{1}{s!(i_1-s)!}
\mbox{ \ \ for \ }1\le s\le i_1.
$$
Using this in (\ref{p2-3}), we obtain that the coefficient of $\hA(g_\bi )$ in $\hH(g)$ is
$\frac{1}{i_1!\cdots i_t!}$. This proves (\ref{p2.1}) and the lemma.
\end{proof}
\begin{remark}\label{re-com}
\rm
Note that a special case of (\ref{p2.1}) is when $g_1=...=g_r$. In
this case we have
\equa{special-case}
{
\hH(g_1,...,g_1)=\frac{1}{r!}\hA(g_1,...,g_1)
\mbox{ \ \ \ (cf.~definition (\ref{denote-M}))}.
%\vs{-6pt}
}
%\vs{-6pt}\PED
\end{remark}
\begin{remark}\label{re-com'}
\rm Using notation (\ref{define-cC}),
formula (\ref{(17)}) can be written as
\equa{for-2}
{
\ch V(\L)=
\mbox{$\sum\limits_{\si\in S^\L}$}
\chiP{\HC{(\si\cdot\L)\lm}}.
}
We also have
\equa{chi-l}
{
\chiP{\NC\l}=
\mbox{$\sum\limits_{\si\in \Sr_r}$}
\chiP{\HC{(\si\cdot\l)\lm}}
\mbox{ \ \ \ for \ }\l\in\NC\L.
}
To prove (\ref{chi-l}),
note that the derivation of (\ref{for-2}) from formula  (\ref{formula})
does not depend on how $\dd^\l_{s,t}$'s are defined. Thus if
we  simply regard $\dd^\l_{s,t}$ as zero for all $s,t$,
then $\TC\l$ coincides
with $\NC\l$, and $S^\l$ becomes $\Sr_r$. Hence (\ref{chi-l})
can be regarded as a special
case of (\ref{for-2}).
Formula (\ref{lemm7-2-1}) can be re-written as
\equa{re-w}
{
\chi(\HC\l)=\mbox{$\sum\limits_{\tau \in C_r}$}\frac{1}{r!}
\left({\,}^{^{\dis r}}_{_{\dis \tau}}{\ssc\,}\right)
(-1)^{\ell(\tau )}\chi(\NC{(\tau\cdot\l)\lm}).
}
Thus formula (\ref{re-w}) is the inverse formula of (\ref{chi-l}).
\PED\end{remark}
\subsection{Kac-Weyl type formula}
While the character formula (\ref{ch-v}) or (\ref{formula}) is extremely useful for understanding structural features
of irreducible $\fg$-modules, such as their resolutions in terms of Kac-modules, it is
not easy to use for purposes
like determining the dimensions of irreducibles. For such purposes, a Kac-Weyl type formula is more desirable.
We now derive such a formula.

For any weight $\l\in \PP $ and any subset $\G\subset\D_1^+$, we define
\equa{define-formula-type}
{
\chi^{\rm BL}_\G(\l)=L_0^{-1}\mbox{$\sum\limits_{w\in W}$}\es(w)e^{w(\l+\rho_0)}
\mbox{$\prod\limits_{\b\in\D_1^+\bs\G}$}(1+e^{-w(\b)}),
}
(which was referred to as the {\it Bernstein-Leites type character} in \cite{VHKT}).
Here and below $L_0,L_1$ are defined in (\ref{l-0}).
For any $\l\in\NC\L$
(not necessarily regular), by (\ref{define-cC}) and
the definition of the Kac-character $\CHIK(\L)$ in (\ref{kac-ch}),
we have
\begin{eqnarray}
\label{fi1}
\nonumber
\chiP{\NC\l}&\!\!\!=\!\!\!&(-1)^{|\L-\l|}
\frac{L_1}{L_0}\ \mbox{$\sum\limits_{w\in W}$}\ \,
\mbox{$\sum\limits_{0\le i_1,...,i_r<\infty}$}\,
(-1)^{\sum_{s=1}^ri_s}
e^{\l+\rho-\sum_{s=1}^ri_s\g_s}
\\
%\nonumber
&\!\!\!=\!\!\!&(-1)^{|\L-\l|}\,
\frac{L_1}{L_0}\ \mbox{$\sum\limits_{w\in W}$}\,
\frac{e^{\l+\rho}}{\prod_{s=1}^r(1+e^{-\g_r})}
\\\nonumber&\!\!\!=\!\!\!&
(-1)^{|\L-\l|}\chi^{\rm BL}_{\G_\L}(\l),
\end{eqnarray}
where $\G_\L=\{\g_1,...,\g_r\}$ defined in (\ref{G-L})
is the set of atypical roots of $\L$
(cf.~(\ref{aty-roots})). In deriving
the last equality one has made use of
the expression of $L_1$ in (\ref{l-0}).

This together with (\ref{for-2}) and (\ref{lemm7-2-1})
(or (\ref{re-w})) proves
\equa{for1}
{\ch V(\L)
=
\sum_{\si\in S^\L,\,\tau \in C_r}{}
\frac{1}{r!}
\left({\,}^{^{\dis r}}_{_{\dis \tau}}{\ssc\,}\right)
(-1)^{|\L-(\tau\cdot(\si\cdot\L)\lm)\lm|+\ell(\tau )}
\chi^{\rm BL}_{_{\G_\L}}\bigl((\tau\cdot(\si\cdot\L)\lm)\lm\bigr).
}
Thus we obtain
the following theorem.
\begin{theorem}
\label{theo7.2}
The formal character  $\ch V(\L)$ of the finite dimensional
irreducible $\fg$-module $V(\L)$ is given by
\\[8pt]
$\ch V(\L)$\\[-6pt]
\equa{for2}{
{\sc\!\!\!\!\!\!\!\!}
=\!\!\!\!\rb{-2pt}{$\dis\sum_{\si\in S^\L,\,\tau\in C_r}{}$}\!\!\!
\frac{1}{r!}\!
\left({\ssc\,}^{^{\dis r}}_{_{\dis \tau}}\right){\!}
(-1)^{|\L-(\tau\cdot(\si\cdot\L)\lm)\lm|+\ell(\tau)}\!
\frac{1}{L_0}\!
\rb{-2pt}{$\dis\sum_{w\in W}\!{\sc\!}{}$}
\es(w)
w{\sc\!}\biggl(\!
e^{(\tau\cdot(\si\cdot\L)\lm)\lm+\rho_0}
\rb{-2pt}{\mbox{$\!\!\!\!\dis\prod\limits_{\b\in\D_1^+\bs \G_\L}\!\!\!\!$}}
(1{\sc\!}+{\sc\!}e^{-\b})
\!\biggr)
.
}
where notations $S^\L,\ \l\lm,\ C_r$ are defined by
$(\ref{s-f-r}),\ (\ref{l-m})$ and \Notation $\ref{Notation3}$
$($see also $(\ref{G-L}),$ $(\ref{l-0}),$
$(\ref{relative}),$
$(\ref{ell-si})$ and
$(\ref{length})$
for other notations$)$.
\end{theorem}
\begin{remark}\label{Remark1}
\rm
Let $\l\in\L+\sum_{s=1}^r\Z\g_s$.
\def\Re#1{{\rm Reg}_{#1}}
Denote the set of regular lexical weights which are $\le\l$ by:
\equa{l-reg}
{
\Re{\l}
=\{\mu\in\NC\L\,|\,\,\mu\le\l,\mbox{ and }\mu\mbox{ is regular and
lexical}\}.
}
Define
\equa{l-ml}
{
\l\bm
=\max\{\mu\in\L+\dsum{s=1}r\Z\g_s\,\,|\,\,\
\Re{\mu}\ =\ \Re{\l}
\mbox{ \ and $\mu$ is lexical}\}.
}
Thus %From the definition one can prove that
$\l\bm$ is obtained from $\l\lm$ by replacing the atypical entries by  (cf.~(\ref{l-m}))
\equan{l-ml'}
{
\shb{\l\bm}{s}{}=\shb{\l\lm}{t(s)}{}, \mbox{ \ \ where \ }
t(s)=\max\{t\ge s\,|\,\,\dd^{\l\lm}_{p,{\ssc\,}t}=1\mbox{ for }s\le p\le t\},
\eqno(\ref{l-ml})'
}
for $s=1,...,r$ (note that although $\l\lm$ is not necessarily regular,
one can still \mbox{define} $\dd^{\l\lm}_{p,{\ssc\,}t}$ as stated in Remark \ref{reme-in-D-r}).
Then in formula (\ref{(16)}) (hence also in formula (\ref{for2})),
for each $\si\in S^\L$, the
$(\si\cdot\L)\lm$
can be replaced by any lexical
weight $\eta$ with
$(\si\cdot\L)\lm\le\eta\le(\si\cdot\L)\bm$.
In particular, we have another character formula:
\\[8pt]
$\ch V(\L)$\\[2pt]
$\dis
=\!\!\!\!\rb{-2pt}{$\dis\sum_{\si\in S^\L,\,\tau\in C_r}{}$}\!\!\!
\frac{1}{r!}
\left({\ssc\,}^{^{\dis r}}_{_{\dis \tau}}\right){\!}
(-1)^{|\L-(\tau\cdot(\si\cdot\L)\bm)\lm|+\ell(\tau)}\!
\frac{1}{L_0}\!
\rb{-2pt}{$\dis\sum_{w\in W}\!{\sc\!}{}$}
\es(w)
w{\sc\!}\biggl(\!
e^{(\tau\cdot(\si\cdot\L)\bm)\lm+\rho_0}
\rb{-2pt}{\mbox{$\!\!\!\!\dis\prod\limits_{\b\in\D_1^+\bs \G_\L}\!\!\!\!$}}
(1{\sc\!}+{\sc\!}e^{-\b})
\!\biggr).
$\hfill$(\ref{for2})'$
%\\[-10pt]\hs{4ex}\ \PED
\end{remark}
\begin{remark}\label{Remark2}
\rm
The main difference between (\ref{for2}) and (\ref{for2})$'$ lies in that
formula (\ref{for2})$'$ keeps the number of
the distinct weights $(\tau\cdot(\si\cdot\L)\bm)\lm$
to be minimal. For example when $\L=0$, then
$(\tau\cdot(\si\cdot\L)\bm)\lm=0\bm$ for all $\tau\in C_r$ and $\si\in
S^\L=\{1\}$ (See Corollary \ref{coro7.2} below).
But all $(\tau\cdot(\si\cdot\L)\lm)\lm$ for $\si=1,\,\tau\in C_r$ are
distinct (in this case $(\si\cdot\L)\lm=0$).
\PED\end{remark}
\begin{coro}\label{coro7.2}
If $\L$ is totally connected
$($cf.~definition after Theorem $\ref{theo7.1})$, then
\begin{eqnarray}
\label{t-con}
\ch(\L)&\!\!\!=\!\!\!&
\frac{1}{r!L_0}%\frac{1}{L_0}
\sum_{w\in W}
\es(w)w\biggl(
e^{\L\bm+\rho_0}\prod_{\b\in\D_1^+\bs \G_\L}(1+e^{-\b})\biggr)
=\frac{1}{r!}\chi^{\rm BL}_{\G_\L}(\L\bm),
\end{eqnarray}
where $\L\bm$ is defined by the way that
$\L\bm+\RHO$ is obtained
from $\L+\RHO$
by replacing all atypical entries $($cf.~$(\ref{denote-t-f}))$ by the largest one.
In
particular, by taking $\L=0$ we obtain the following denominator formula
\equa{de-for}
{
L_0=\frac{1}{r_0!}\sum_{w\in W}\es(w)w\biggl(e^{0\bm+\rho_0}
\prod_{\b\in\D_1^+\bs\G_0}(1+e^{-\b})\biggr),
}
where $r_0=\min\{m,n\}$,
$\G_0=\{\g^0_1=\es_m-\d_1,...,
\g^0_{r_0}=\es_{m+1-r_0}-\d_{r_0}\}$ and
$$\mbox{$
0\bm=\sum\limits_{s=1}^{r_0}(r_0-s)\g^0_s=
(0,...,0,1,...,r_0-1\,|\,r_0-1,...,1,0,...,0)
\mbox{\ \ $($cf.~notation $(\ref{weight1}))$.}
$}$$
\end{coro}
\begin{proof}
If $\L$ is totally connected, then $S^\L=\{1\}$ and (\ref{t-con}) follows from
(\ref{special-case}) and (\ref{(17)}) with $(\si\cdot\L)\lm$ replaced by
$\L\bm$. Observe that when $\L=0$, we have $\ch V(\L)=1$ and
$\G_\L=\G_0$. Thus we obtain (\ref{de-for}).
\end{proof}
\begin{coro}\label{coro7.2'}
If $\L$ is totally disconnected
$($cf.~definition after Theorem $\ref{theo7.1})$, then
$$
\ch(\L)
=\frac{1}{L_0}
\sum_{w\in W}
\es(w)w\biggl(e^{\L+\rho_0}\prod_{\b\in\D_1^+\bs \G_\L}(1+e^{-\b})\biggr)
=\chi^{\rm BL}_{\G_\L}(\L).
$$
\end{coro}
\begin{proof}
If $\L$ is totally disconnected, then $S^\L=\Sr_r$
and the result follows from the fact that
(\ref{lemm7-2-1}) is the inverse formula of (\ref{chi-l}) (cf.~(\ref{fi1})).
\end{proof}
\subsection{Dimension formula} An important application of Theorem \ref{theo7.2} is the derivation
of a dimension formula.
\begin{theorem}
\label{theo7.3}
The dimension $\dim V(\L)$ of the finite dimensional
irreducible $\fg$-module $V(\L)$ is given by
\begin{eqnarray}
\nonumber
\!\!\!\!&\!\!\!\!\!\!\!\!&\!\!\!\!\!\!\!\!\!\!\!\!\!\!\!\!\!\!
{\rm dim\,} V(\L)\\
\!\!\!\!&\!\!\!\!\!\!\!\!&\!\!\!\!\!\!\!\!\!\!\!\!\!\!\!\!\!\!
=\!\!
\sum_{^{\sc\si\in S^\L,{\ssc\,}\tau \in C_r}_{\sc B\subset\D_1^+\bs\G_\L}}\!
\frac{1}{r!}
\left({\ssc\,}^{^{\dis r}}_{_{\dis \tau}}\right)
(-1)^{|\L-(\tau\cdot(\si\cdot\L)\lm)\lm|+\ell(\tau)}\!
{\!}
\prod_{\a\in\D_0^+}\!\!
\frac{(\a,\rho_0+(\tau\cdot(\si\cdot\L)\lm)\lm-\sum_{\b\in B}\b)}{(\a,\rho_0)}.
\end{eqnarray}
\end{theorem}
\begin{proof}
Regard an element of $\varepsilon$ (cf.~(\ref{epsilon})) as a function on $\fh^*$ such that
the evaluation of $e^\l$ on $\mu$ is
$e^\l(\mu)=e^{(\l,\mu)}$ for $\mu\in\fh^*$. Then
$\ch V(\L)\in \varepsilon$ is a function on $\fh^*$, and
\equa{dim-V}
{
\dim V(\L)=\lim_{x\to0}\,(\ch V(\L))(x\rho_0).
}
First we calculate
%the evaluation of $\chi^{\rm BL}_\G(\l)$ on $x\rho_0$
$\lim_{x\to0}\chi^{\rm BL}_\G(\l)(x\rho_0)$ (cf.~(\ref{define-formula-type})).
Using
$$\mbox{$
\prod\limits_{\b\in\D_1^+\bs\G_\L}(1+e^{-\b})=
\sum\limits_{B\subset\D_1^+\bs\G_\L}e^{-\sum_{\b\in B}\b},
$}$$
and
the well-known denominator formula
\equa{de-f}{
L_0=\mbox{$\sum\limits_{w\in W}$}\es(w)e^{w(\rho_0)}
\mbox{ \ \ \ (cf.~(\ref{l-0})),}
}
we have
\begin{eqnarray*}
\!\!\lim_{x\to0}
\chi^{\rm BL}_\G(\l)(x\rho_0)
&\!\!=\!\!&
\lim_{x\to0}\,
\mbox{$\sum\limits
_{B\subset\D_1^+\bs\G_\L,\,w\in W}$}\es(w)
\frac{e^{(w(\l+\rho_0-\sum_{\b\in B}\b),x\rho_0)}}{L_0(x\rho_0)}
\\
&\!\!=\!\!&
\lim_{x\to0}\,
\mbox{$\sum\limits_{B\subset\D_1^+\bs\G_\L}$}
\frac{L_0(x(\l+\rho_0-\sum_{\b\in B}\b))}{L_0(x\rho_0)}
\\%\end{eqnarray*}\begin{eqnarray*}
&\!\!=\!\!&
\lim_{x\to0}\,
\mbox{$\sum\limits_{B\subset\D_1^+\bs\G_\L}\ \prod\limits_{\a\in\D_0^+}$}
\frac{e^{(\a/2,x(\l+\rho_0-\sum_{\b\in B}\b))}-e^{(-\a/2,x(\l+\rho_0
-\sum_{\b\in B}\b))}}{e^{(\a/2,x\rho_0)}-e^{(-\a/2,x\rho_0)}}
\\%\end{eqnarray*}\begin{eqnarray*}
&\!\!=\!\!&
\mbox{$\sum\limits_{B\subset\D_1^+\bs\G_\L}\ \prod\limits_{\a\in\D_0^+}$}
\frac{(\a,\l+\rho_0-\sum_{\b\in B}\b)}{(\a,\rho_0)},
\end{eqnarray*}
where the second equality follows from (\ref{de-f}) and the fact that
$(w(\l),\mu)=(\l,w(\mu))$ for $w\in W,\,\l,\mu\in\fh^*$.
This together with (\ref{dim-V}) and (\ref{for1}) gives the result.
\end{proof}
\begin{remark}
\rm
As far as we are aware, there exists no dimension formula for the finite dimensional irreducible
$\fg$-modules in the literature
except the one for singly atypical modules given by van der Jeugt \cite{V1}.
%
%
%{\bf Note: two ref \cite{VHKT0} and \cite{V1} added.}
%
%
\PED\end{remark}

%%%%%%%%%%%%%%%%%%%%%%%%%%%%%%%%%%%%%%%
{\small\vskip4pt \noindent{\bf Acknowledgements}.
We wish to thank Professor V. Kac for
suggesting important improvements to the paper.
Financial support from the Australian Research Council is
gratefully acknowledged. Su is also supported by NSF grant 10171064 of
China, EYTP and TCTPFT grants of Ministry of Education of China.}

%%%%%%%%%%%%%%%%%%%%%%%%%%%%%%%%%%%%%%%%%%%%%%%%%%%%%%%%


\begin{thebibliography}{9999}
%%%%%%%%%%%%%%%%%%%%%%%%%%%%%%%%%%%%%%%%%%%%%%%%%%%%%%%%
\bibitem{BL} I.N.~Bernstein, D.A.~Leites,
{\em A formula for the characters of the irreducible finite-dimensional
representations of Lie superalgebras of series $gl$ and $sl$},
C.~R.~Acad.~Bulgare Sci.~{\bf33} (1980), 1049-1051.

\bibitem{B} J.~Brundan, {\em Kazhdan-Lusztig polynomials and
   character formulae for the Lie superalgebra $gl(m|n)$},
   J.~Amer.~Math.~Soc.~{\bf 16} (2002) 185-231.

\bibitem{B1} J.~Brundan, {\em Kazhdan-Lusztig polynomials and character
formulae for the Lie superalgebra $q(n)$},  Adv.~Math. {\bf 182}  (2004) 28-77.


\bibitem{CZ} S.-J. Cheng and R.B. Zhang, {\em  An analogue of
Kostant's ${\mathfrak u}$-cohomology
formula for the general linear superalgebra}.
International Mathematics Research Notices {\bf1} (2004) 31-53.

\bibitem{HKV} J.W.B.~Hughes, R.C.~King and J.~van der Jeugt,
{\em On the composition factors of Kac modules for the Lie superalgebras $sl(m|n)$},
J.~Math.~Phys.~{\bf33} (1992) 470-491.

\bibitem{Kac0} V.G.~Kac, {\em Classification of simple Lie superalgebras},
Funct.~Anal.~Appl. {\bf9} (1975) 263-265.

\bibitem{Kac} V.G.~Kac, {\em Lie superalgebras},
   Adv.~Math.~{\bf 26} (1977) 8-96.

\bibitem{Kac1} V.G.~Kac, {\em Characters of typical representations of
  classical Lie superalgebras}, Comm.~Alg.~{\bf5} (1977) 889-897

\bibitem{Kac2} V.G.~Kac, {\em Representations of classical Lie superalgebras},
  Lect. Notes Math. {\bf676} (1978) 597-626.

\bibitem{PS97} I.B.~Penkov and V.~Serganova, {\em Characters of
irreducible $G$-modules
   and cohomology of $G/P$ for the Lie supergroup $G=Q(N)$},
   Algebraic geometry, 7, J.~Math.~Sci.~(New York) {\bf 84} (1997), no.~5,
   1382-1412.


\bibitem{Se96} V.~Serganova, {\em Kazhdan-Lusztig polynomials and
       character formula for the Lie superalgebra $gl(m|n)$},
       Selecta Math.~{\bf 2} (1996) 607-654.

\bibitem{Se98} V.~Serganova,
  {\em Characters of irreducible representations of simple Lie superalgebras},
  Proceedings of the International Congress of Mathematicians 1998, Berlin,
  Vol.~II, Documenta Mathematica, Journal der Deutschen
  Mathematiker-Vereinigung, pp.~583-593.

\bibitem{V} J.~van der Jeugt, {\em Character formulae for the
  Lie superalgebra $C(n)$}, Comm.~Algebra {\bf19} (1991) 199-222.

\bibitem{V1} J.~van der Jeugt, {\em Dimension formulas for
the Lie superalgebras $sl(m|n)$}, J.~Math.~Phys.~{\bf 36} (1995) 605-611.


\bibitem{VHKT0} J.~van der Jeugt, J.W.B.~Hughes, R.C.~King and J.~Thierry-Mieg,
{\em A character formula for singly atypical modules of the Lie superalgebra
$sl(m/n)$}, Comm.~Alg.~{\bf 19} (1991) 199-222.



\bibitem{VHKT} J.~van der Jeugt, J.W.B.~Hughes, R.C.~King and
J.~Thierry-Mieg, {\em Character formulas for irreducible modules of
the Lie superalgebras $sl(m|n)$}, J.~Math.~Phys.~{\bf31} (1990) 2278-2304.

\bibitem{VZ} J.~van der Jeugt and R.B.~Zhang, {\em Characters and composition
 factor multiplicities for the Lie superalgebra $gl(m|n)$},
  Lett.~Math.~Phys.~{\bf 47} (1999) 49-61.

\bibitem{Zou} Y.M.~Zou, {\em Categories of finite-dimensional weight
modules over type I classical Lie superalgebras},  J.~Algebra {\bf 180} (1996) 459-482.

\end{thebibliography}
\end{document}